\documentclass[letterpaper, 10pt,psfig, conference]{ieeeconf}

\usepackage{amssymb}
\usepackage{graphicx}
\usepackage{amsmath}%
\usepackage{amsfonts}
\usepackage{epsfig}

\newcommand{\nc}{\newcommand}
\nc{\dd}{d_{\Delta}}
\nc{\uk}{\underline{k}}
\nc{\bk}{\bar{k}}

\nc{\that}{\hat{\theta}}

\nc{\rhod}{{\rho_{\delta}}}
\nc{\rhot}{{\rho_{\delta} ( \phi (t) , e(t+1))}}
\nc{\rhoj}{{\rho_{\delta} ( \phi (j) , e(j+1))}}
\nc{\rhotm}{{\rho_{\delta} ( \phi (t) , e(t+1))}}
\nc{\rhojm}{{\rho_{\delta} ( \phi (j) , e(j+1))}}

\nc{\bphi}{\bar{\phi}}
\nc{\tphi}{\tilde{\phi}}
\nc{\tct}{{\tilde{\check{\theta}}}}
\nc{\vc}{\check{V}}
\nc{\ttt}{\tilde{\theta}}
\nc{\D}{{{\bf D}}}

\nc{\xnun}{{\tbo{x(0)}{\hat{u}^0 (0)}}}
\nc{\xkuk}{{\tbo{x(kT)}{\hat{u}^0(kT)}}}

\nc{\blz}{\bar{\lambda}_0}
\nc{\bgz}{\bar{\gamma}_0}
\nc{\bh}{\bar{h}}
\nc{\bn}{\bar{n}}
\nc{\barki}{{\bar{k}_i}}
\nc{\barkio}{{\bar{k}_{i-1}}}

\nc{\bub}{{\bar{\underline{B}}}}
\nc{\bcub}{{\bar{\underline{C}}}}

\nc{\bb}{{\bar{B}}}
\nc{\bc}{{\bar{C}}}
\nc{\ub}{{\underline{B}}}
\nc{\uc}{{\underline{C}}}
\nc{\ucb}{{\underline{\cal B}}}
\nc{\ucc}{{\underline{\cal C}}}
\nc{\htxe}{{\hat{\tilde{x}}^{\eps}}}
\nc{\htue}{{\hat{\tilde{u}}^{\eps}}}
\nc{\htye}{{\hat{\tilde{y}}^{\eps}}}
\nc{\rmr}{{\R^{\tilde{m} \times \tilde{r}}}}
\nc{\kde}{{ \tilde{K}_{\delta}^{\eps}}}
\nc{\kdep}{{ \tilde{K}_{\delta}^{\eps}(p)}}
\nc{\omc}{{{\cal O}_{m}(C,A)}}
\nc{\tkt}{{ [kT , (k+1)T)}}
\nc{\mbl}{{\frac{\bar{m}}{\ell}}}
\nc{\calo}{{\cal O}}
\nc{\calc}{{\cal C}}
\nc{\acl}{{A_{cl}}}
\nc{\acla}{{A_{cl(A,B,C)}}}
\nc{\aclep}{\tilde{A}_{cl}^{\eps} (p)}
\nc{\gabc}{{\gamma_{(A,B,C)}}}
\nc{\gabcb}{{\gamma_{( \bar{A} , \bar{B} , \bar{C} )}}}
\nc{\gabcbt}{{\tilde{\gamma}_{(\bar{A} , \bar{B} , \bar{C})}}}
\nc{\dabc}{{\delta_{(A,B,C)}}}
\nc{\dabcb}{{\delta_{( \bar{A} , \bar{B} , \bar{C} )}}}
\nc{\labc}{{\lambda_{(A,B,C)}}}
\nc{\labcb}{{\lambda_{(\bar{A}, \bar{B}, \bar{C})}}}
\nc{\acleb}{{A_{cl( \bar{A} , \bar{B} , \bar{C} )}^{\eps}}}
\nc{\abcp}{{(A,B,C) \in {\cal P}}}
\nc{\abcbp}{{(\bar{A}, \bar{B}, \bar{C}) \in {\cal P}}}
\nc{\hxe}{{\hat{x}^{\eps}}}
\nc{\hue}{{\hat{u}^{\eps}}}

\nc{\eabcb}{{\varepsilon_{(\bar{A}, \bar{B}, \bar{C})}}}
\nc{\feabc}{{F_{(A,B,C)}^{\varepsilon}}}
\nc{\fabce}{{F_{(A,B,C)}^{\varepsilon}}}
\nc{\fabcb}{{F_{(\bar{A}, \bar{B}, \bar{C})}}}
\nc{\fabcbe}{{F_{(\bar{A}, \bar{B}, \bar{C})}^{\varepsilon}}}
\nc{\fabc}{{F_{(A,B,C)}}}
\nc{\oh}{{\cal O}}
\nc{\xbb}{{\bar{\bar{x}}}}
\nc{\bp}{{\bar{p}}}
\nc{\fguy}{{\cal F}}
\nc{\eps}{{\varepsilon}}
\nc{\feps}{{\hat{f}_{\eps}}}
\nc{\xbeps}{{\bar{x}^{\eps}}}
\nc{\xheps}{{\hat{\bar{x}}^{\eps}}}
\nc{\xhe}{{\hat{{x}}^{\eps}}}
\nc{\xte}{{\tilde{{x}}^{\eps}}}
\nc{\uhe}{{\hat{{u}}^{\eps}}}
\nc{\ute}{{\tilde{{u}}^{\eps}}}

\nc{\nabcb}{{{\cal N}_{(\bar{A}, \bar{B}, \bar{C})}}}
\nc{\nabcbt}{{\tilde{\cal N}_{(\bar{A}, \bar{B}, \bar{C})}}}

\nc{\bl}{{\bar{\lambda}}}
\nc{\R}{{\bf R}}
\nc{\C}{{\bf C}}
\nc{\Z}{{\bf Z}}
\nc{\N}{{\bf N}}
\nc{\linf}{{l_{\infty}}}
\nc{\Fb}{{{F} }}
\nc{\tb}{{\bar{t}}}
\nc{\tub}{{\underline{t}}}
\nc{\umguy}{{\| u_m \|_{\infty}}}

\nc{\stable}   {\mbox{SJVL-stable}}
\nc{\Rica}     {{R\overline{ic}}}
\nc{\Ricb}     {\overline{R\overline{ic}}}
\nc{\B}        {{\cal B} }
\nc{\Ltplusd}  {{\cal L}_2  [0, \infty ) }
\nc{\Ltplus}   {{\cal L}_{2+} }
\nc{\Pplus}    {P_{+} }
\nc{\Ltminusd} {{\cal L}_2 (-\infty ,0] }
\nc{\Ltminus}  {{\cal L}_{2-} }
\nc{\Pminus}   {P_{-} }
\nc{\Hf}       {{\cal H}_{\infty} }
\nc{\Ht}       {{\cal H}_{2} }
\nc{\Htp}      {{\cal H}_{2}^{\perp} }
\nc{\RHt}      {{\cal RH}_{2} }
\nc{\RHtp}     {{\cal RH}_{2}^{\perp} }
\nc{\RHf}      {{\cal RH}_{\infty} }
\nc{\Lf}       {{\cal L}_{\infty} }
\nc{\Lt}       {{\cal L}_2 }
\nc{\RLf}      {{\cal RL}_{\infty} }
\nc{\Xminus}   {{\cal X}_- }
\nc{\Xplus}    {{\cal X}_+ }
\nc{\be}       {\begin{equation}}
\nc{\ee}       {\end{equation}}
\nc{\bea}      {\begin{eqnarray}}
\nc{\eea}      {\end{eqnarray}}
\nc{\ba}       {\begin{array}}
\nc{\ea}       {\end{array}}
\nc{\p}        { ^\prime }
\nc{\Xf}       { X_{\infty} }
\nc{\Yf}       { Y_{\infty} }
\nc{\Ymf}      { Y_{tmp } }
\nc{\Hamf}     { H_{\infty} }
\nc{\Sf}       { S_{\infty}}         % PAI 061089
\nc{\Tf}       { T_{\infty}}         % PAI 061089
\nc{\Ff}       { F }
\nc{\Jf}       { J_{\infty} }
\nc{\Jmf}      { J_{tmp } }
\nc{\Ax}       { A_{X} }
\nc{\AFf}      { A_{F} }
\nc{\AFt}      { A_{F_0} }
\nc{\ALt}      { A_{L_2} }
\nc{\CoFf}     { C_{1F} }
\nc{\CoFt}     { C_{1F_0} }
\nc{\BoLt}     { B_{1L_2} }
\nc{\CFt}      { C_{F_2} }
\nc{\sU}       { {\cal U} }
\nc{\Dp}       {D_{\perp}}
\nc{\tR}       {\tilde{R}}
\nc{\tD}       {\tilde{D}}
\nc{\tDp}      {\tD_{\perp}}
\nc{\Doo}      {D_{11}}
\nc{\Dto}      {D_{21}}
\nc{\Dbo}      {D_{\bullet1}}
\nc{\Dob}      {D_{1\bullet}}
\nc{\Fl}       {{\cal F}_{\ell}}
\nc{\fFl}      {\Fl(G,K)}
\nc{\fFlf}     {\|\fFl\|_{\infty}}
\nc{\Gtmp}     {G_{\rm tmp}}
\nc{\Jtmp}     {J_{\rm tmp}}
\nc{\Xtmp}     {X_{\rm tmp}}
\nc{\Ytmp}     {Y_{\rm tmp}}
\nc{\Ftmp}     {F_{\rm tmp}}
\nc{\Ffo}      {F_1}
\nc{\Fft}      {F_2}
\nc{\Fo}       {F^o}
\nc{\Lo}       {L^o}
\nc{\Tot}      {(T_1^{\prime}T_1)^{-1}}
\nc{\Bottil}   {(B_{12}-B_{22}D_{12}^{\prime}D_{11})}
\nc{\Bootil}   {(B_{12}-B_{22}D_{12}^{\prime}D_{11})}
\nc{\Bottilp}  {(B_{12}^{\prime}-\Doo^{\prime}\Dot B_{22}^\prime)}
\nc{\Bootilp}  {(B_{11}^{\prime}-\Doo^{\prime}\Dot B_{21}^\prime)}
\nc{\DoDpDp}   {\Doo^{\prime}\Dp\Dp^{\prime}}
\nc{\xh}       {\hat{x}}
\nc{\xhd}      {\dot{\hat{x}}}
\nc{\tB}       {\tilde{B}}
\nc{\eqdent}   {\hspace{6ex}}
\nc{\tallarray}{\renc{\arraystretch}{1.15}}
\nc{\normarray}{\renc{\arraystretch}{1.0}}

\nc{\tbth}[6]{
  \left[ \begin{array}{ccc}
       #1 & #2 & #3 \\ #4 & #5 & #6 
       \end{array} \right] }
\nc{\thbt}[6]{
  \left[ \begin{array}{cc}
       #1 & #2  \\ #3 & #4 \\ #5 & #6 
       \end{array} \right] }
\nc{\thbth}[9]{
  \left[ \begin{array}{ccc}
       #1 & #2 & #3 \\ #4 & #5 & #6 \\ #7 & #8 & #9
       \end{array} \right] }
\nc{\tbt}[4]{
  \left[ \begin{array}{cc}
       #1 & #2 \\ #3 & #4
       \end{array} \right] }
\nc{\tbo}[2]{
  \left[ \begin{array}{c}
       #1 \\ #2
       \end{array} \right] }
\nc{\tbf}[8]{
  \left[ \begin{array}{cccc}
       #1 & #2 & #3 & #4 \\
       #5 & #6 & #7 & #8
       \end{array} \right] }
\nc{\thbo}[3]{
  \left[ \begin{array}{c}
       #1 \\ #2 \\ #3
       \end{array} \right] }
\nc{\obf}[4]{
  \left[ \begin{array}{cccc}
       #1 &  #2 &  #3 &  #4
       \end{array} \right] }
\nc{\obfi}[5]{
  \left[ \begin{array}{ccccc}
       #1 &  #2 &  #3 &  #4 & #5
       \end{array} \right] }
\nc{\obs}[6]{
  \left[ \begin{array}{cccccc}
       #1 &  #2 &  #3 &  #4 & #5 & #6
       \end{array} \right] }
\nc{\sbo}[6]{
  \left[ \begin{array}{c}
       #1 \\  #2 \\  #3 \\  #4 \\ #5 \\ #6
       \end{array} \right] }
\nc{\fbo}[4]{
  \left[ \begin{array}{c}
       #1 \\  #2 \\  #3 \\  #4
       \end{array} \right] }
\nc{\fibo}[5]{
  \left[ \begin{array}{c}
       #1 \\  #2 \\  #3 \\  #4 \\ #5
       \end{array} \right] }
\nc{\obt}[2]{
  \left[ \begin{array}{cc}
       #1 & #2
       \end{array} \right] }
\nc{\obth}[3]{
  \left[ \begin{array}{ccc}
       #1 & #2 & #3
       \end{array} \right] }
\nc{\tfmat}[4]{
  \left[ \begin{array}{c|c}
       #1 & #2 \\ \hline
       #3 & #4
       \end{array} \right] }
\nc{\LFT}[6]{
  \begin{picture}(60,42)
    \put(20,25){\framebox(20,16){#1}} % P
    \put(24,5){\framebox(12,10){#2}}  % K
    \put(10,40){\makebox(0,0){#3}}  % z
    \put(6,20){\makebox(0,0){#4}}   % y
    \put(50,40){\makebox(0,0){#5}}  % w
    \put(54,20){\makebox(0,0){#6}}  % u
    \put(20,37){\vector(-1,0){20}}
    \put(60,37){\vector(-1,0){20}}
    \put(50,29){\vector(-1,0){10}}
    \put(10,29){\line(1,0){10}}
    \put(36,10){\line(1,0){14}}
    \put(10,10){\vector(1,0){14}}
    \put(10,10){\line(0,1){19}}
    \put(50,10){\line(0,1){19}}
  \end{picture}}
\nc{\lft}[6]{
  \begin{picture}(30,21)
    \put(10,11){\framebox(10,10){#1}}   % P
    \put(11,0){\framebox(8,6){#2}}              % K
    \put(5,21){\makebox(0,0){#3}}               % z
    \put(5,8){\makebox(0,0)[r]{#4}}             % y
    \put(25,21){\makebox(0,0){#5}}              % w
    \put(25,8){\makebox(0,0)[l]{#6}}    % u
    \put(10,19){\vector(-1,0){10}}
    \put(30,19){\vector(-1,0){10}}
    \put(24,13){\vector(-1,0){4}}
    \put(6,13){\line(1,0){4}}
    \put(19,3){\line(1,0){5}}
    \put(6,3){\vector(1,0){5}}
    \put(6,3){\line(0,1){10}}
    \put(24,3){\line(0,1){10}}
  \end{picture}}
\nc{\LFTcenter}[6]{
  \begin{center}
    \LFT{#1}{#2}{#3}{#4}{#5}{#6}
  \end{center}}
\nc{\lftcenter}[6]{
  \begin{center}
    \lft{#1}{#2}{#3}{#4}{#5}{#6}
  \end{center}}
\nc{\lfte}[7]{
  \begin{picture}(80,21)
    \put(0,0){\lft{#1}{#2}{#3}{#4}{#5}{#6}}
    \put(50,10){#7}
  \end{picture}}
\nc{\lftecrunch}[7]{
  \begin{picture}(80,21)
    \put(0,0){\lft{#1}{#2}{#3}{#4}{#5}{#6}}
    \put(40,10){#7}
  \end{picture}}
\nc{\lftss}[6]{
  \begin{picture}(64,42)
    \put(14,25){\framebox(32,16){#1}}
    \put(24,5){\framebox(12,10){#2}}
    \put(8,40){\makebox(0,0){#3}}
    \put(6,20){\makebox(0,0){#4}}
    \put(52,40){\makebox(0,0){#5}}
    \put(54,20){\makebox(0,0){#6}}
    \put(14,37){\vector(-1,0){14}}
    \put(60,37){\vector(-1,0){14}}
    \put(50,29){\vector(-1,0){4}}
    \put(10,29){\line(1,0){4}}
    \put(36,10){\line(1,0){14}}
    \put(10,10){\vector(1,0){14}}
    \put(10,10){\line(0,1){19}}
    \put(50,10){\line(0,1){19}}
  \end{picture}}
\nc{\lftsscenter}[6]{
  \begin{center}
    \lftss{#1}{#2}{#3}{#4}{#5}{#6}
  \end{center}}

\newtheorem{theorem}    {Theorem}
\newtheorem{prop}{Proposition}

\newtheorem{lemma}      {Lemma}

\newtheorem{remark}     {Remark}

\begin{document}
\title{\Large \bf {Classical Pole Placement Adaptive Control Revisited: 
Exponential Stabilization}
\author{
Daniel E. Miller\footnote{This research was supported
by a grant from the Natural Sciences Research Council of Canada.} \\
Dept.\ of Elect.\ and Comp.\ Eng.  \\
University of Waterloo, Waterloo, ON  \\
Canada \ N2L 3G1 \\
(miller@uwaterloo.ca) 
}}

\thispagestyle{empty}
\date{\today}
\maketitle
\thispagestyle{empty}

\begin{abstract}

While the original classical parameter adaptive controllers
did not handle noise or unmodelled dynamics well, redesigned versions
were proven to have some tolerance; however, exponential stabilization
and a bounded gain on the noise
was rarely proven.
Here we consider a classical pole placement adaptive controller using the 
original projection algorithm rather than the commonly modified version; we impose
the
assumption that the plant parameters lie in a convex,
compact set.
% and that the parameter estimates are projected onto that set at ever step.
We demonstrate that the closed-loop system exhibits very desireable closed-loop behaviour:
there are linear-like convolution bounds on the closed
loop behaviour, which confers exponential stability and a bounded
noise gain, and can be leveraged to prove tolerance to unmodelled dynamics and
plant parameter variation.
We emphasize that there is no persistent excitation requirement of any sort;
the improved performance arises from the vigilant nature of the
parameter estimator.

\end{abstract}


\vspace{0.3cm}

\section{Introduction}

Adaptive control is an approach used to deal with systems with uncertain
or time-varying parameters.
The classical adaptive controller consists of a linear time-invariant (LTI)
compensator together with a tuning mechanism to adjust the compensator
parameters to match the plant.
The first general proofs that adaptive controllers could work
came around 1980, e.g. see \cite{morse1978},
\cite{morse1980}, \cite{goodwin1980}, \cite{Narendra1980}, and
\cite{Narendra1980_pt2}. However,
such controllers were typically not robust to unmodelled dynamics, did not tolerate
time-variations well, 
and did not handle noise or disturbances well, e.g. see \cite{rohrs}.
During the following two decades a great deal of effort was made to
address these shortcomings. The most common approach was to make
small controller design changes,
such as the use of 
signal normalization, deadzones, and $\sigma-$modification,
to ameliorate these issues,
%so that the resulting controllers tolerated a degree of
%unmodelled dynamics, slow time-variations, and noise, 
e.g. see \cite{rick}, \cite{rick2},
\cite{Tsakalis4}, \cite{kreiss}, \cite{Ioa86}.
Indeed, simply using projection
(onto a convex set of admissible parameters)
has proved quite powerful,
and
the resulting controllers typically provide a bounded-noise bounded-state
property, as well as tolerance of some degree of unmodelled dynamics and/or time-variations,
e.g. see \cite{ydstie}, \cite{ydstie2}, \cite{naik}, \cite{wenhill}, \cite{wen} and \cite{hanfu}.
Of course, it is clearly desireable that the closed-loop system exhibit LTI-like
system properties, such as a bounded gain and exponential stability.
As far as the author is aware,
in the classical approach to adaptive control
a bounded gain on the noise 
\footnote{Since the closed-loop system is nonlinear, a bounded-noise bounded-state property
does not automatically imply a bounded gain on the noise.}
is proven only in \cite{ydstie2};
however, 
%a convolution-like bound on the behaviour is not proven (this is
%critical to proving that the approach can tolerate slow time-variations - see Section 6), 
a crisp exponential bound on the effect of the initial condition
is not provided, and a 
minimum phase assumption is imposed.
While it is possible to prove a form of exponential stability if the
reference input is sufficiently persistently exciting, e.g. see
\cite{Nar87}, this places a stringent requirement on an exogenous input.

There are several non-classical approaches to adaptive control which provide
LTI-like system properties. First of all,
in \cite{barmish} and \cite{miller89}
a logic-based switching approach was used to switch between a predefined list
of candidate controllers;
while exponential stability is proven, the transient behaviour
can be quite poor and a bounded gain on the noise is not proven.
A more sophisticated logic-based approach,
labelled
Supervisory Control, was proposed by Morse; here
a supervisor switches in an efficient way between candidate controllers - see
\cite{morse96}, \cite{morse97},
\cite{liberzon1}, \cite{vu2} and \cite{liberzon2}.
In certain circumstances a bounded gain on the noise 
can be proven - see \cite{vu} and 
the Concluding Remarks section of \cite{morse97}.
A related approach, called localization-based switching adaptive control,
uses a falsification approach to prove exponential stability as well
as a degree of tolerance of disturbances, e.g. see \cite{rick3}.

Another non-classical 
approach, proposed by the author, is based on periodic
estimation and control: rather than estimate the plant or controller
parameters, the goal is to estimate what the control signal 
would be if the
plant parameters and plant state were known and the `optimal controller'
were applied.
%The
%goal is to develop controllers for which there are linear-like bounds on the closed loop behaviour:
Exponential stability and a bounded gain on the noise is achieved, as well as
near optimal performance, e.g. see
\cite{miller03}, \cite{miller06}, and \cite{julie11};
a degree of unmodelled dynamics
and time variations can be allowed.
The cost of these desireable features is that the noise gain
increases dramatically the closer that one gets to optimality.

In this paper we consider the discrete-time setting and we propose an alternative approach to obtaining LTI-like system properties.
%Indeed, we are able to obtain linear-like convolution bounds
%on the closed-loop behaviour, which automatically guarantees exponential stability and a bounded
%gain on the noise but also can be leveraged to
%prove a degree of tolerance to
%parameter variation as well as unmodelled dynamics.
We return to a common approach in classical adaptive control - the use of the projection
algorithm together with the Certainty Equivalence Principle.
In the literature it is the norm to use a modified version of the ideal Projection
Algorithm in order to avoid division by zero;
\footnote{
An exception is the work of Ydstie \cite{ydstie}, \cite{ydstie2},
who considers the ideal Projection Algorithm as a special case; however, 
a crisp bound on the effect of the initial condition is not proven
and
a minimum phase assumption is imposed.}
it turns out that
an unexpected consequence of this minor adjustment is that some inherent
properties of the scheme are destroyed.
Here we use the original version of the Projection Algorithm coupled with a pole placement
Certainty Equivalence based controller. We 
obtain linear-like convolution bounds
on the closed-loop behaviour, which immediately confers
exponential stability and a bounded gain on the noise; such convolution
bounds are, as far as the author is aware, a first in adaptive control, and it
allows us to use a modular approach to
analyse robustness and tolerance to time-varying 
parameters.
%We can use all of the intuition that we have developed
%for LTI systems in the adaptive setting, something which
%has not been present using other techniques. 
To this end, the results will be presented in
a very pedagogically desireable fashion: we first deal with the ideal plant
(with disturbances);
we then leverage that result to prove that 
a large degree of time-variations is tolerated;
%commonly considered in the adaptive control literature; 
we then demonstrate that the approach tolerates a degree of unmodelled
dynamics, in a way familiar to those versed in the analysis of LTI
systems.
%
%This paper provides a new approach to dealing with classical
%discrete-time adaptive control, with new proof and analysis techniques.
%We hope that it will lead to a renewal of interest in the field.

In a recent short paper we consider the first order case \cite{scl}.
Here we consider the general case, which requires much more sophisticated analysis and proofs.
Furthermore, in comparison to \cite{scl}, here we (i) present a more
general estimation algorithm, which alleviates  
the classical concern about dividing by zero, (ii)
prove that the controller achieves the objective
in the presence of a more general class of time-variations, and
(iii) prove robustness to unmodelled dynamics.
An early version of this paper has been submitted to a conference \cite{CCTA}.

Before proceeding we present some mathematical preliminaries.
Let $\Z$ denote the set of integers,
$\Z^+$ the set of non-negative integers, $\N$
the set of natural numbers, $\R$ the set of real numbers, and
$\R^+$ the set of
non-negative real numbers.
We let $\D^0$ denote the open unit disk of the complex plane.
We use the Euclidean $2$-norm for vectors and the corresponding induced
norm for matrices, and denote the norm of a vector or matrix
by $\| \cdot \|$.
We let $\linf ( \R^n)$ denote the set of $\R^n$-valued bounded sequences;
we define the norm of $u \in \linf ( \R^n)$
by $\| u \|_{\infty} := \sup_{k \in \Z} \| u(k) \|$.
Occasionally we will deal with a map $F: \linf ( \R^n) \rightarrow
\linf ( \R^n)$; the gain is given by
$\sup_{ u \neq 0 } \frac{\| Fu \|_{\infty}}{\| u \|_{\infty} }$
and denoted by $\| F \|$.
With $T \in \Z$, the truncation operator $P_T: \linf ( \R^n) \rightarrow
\linf ( \R^n)$ is defined by
\[
(P_T x) (t) = \left\{ \begin{array}{ll}
x(t) & \;\; {t  \leq T} \\
0  & \;\; t > T.
\end{array}
\right .
\]
We say that the map $F: \linf ( \R^n) \rightarrow
\linf ( \R^n)$ is causal if $P_T F P_T = P_T F$ for every $T \in \Z$.
%and strictly causal if $P_T F P_{T-1} = P_T F$ for every $T \in \Z$.

If ${\cal S} \subset \R^p$ is a convex and compact set, we
define $\| {\cal S} \| := \max_{x \in {\cal S} } \| x \|$ and the
function $\pi_{\cal S} : \R^p \rightarrow
{\cal S}$ denotes the projection onto ${\cal S}$; it is well-known that
$\pi_{\cal S}$
is well-defined.
%With $\eps > 0$ and $c>0$, we let $s ( {\cal S} , c_0 , \eps )$ denote the set of sequences in
%$x \in \linf ( \R^n) $
%taking values in ${\cal S}$ and
%satisfying, for every $t_2 \in \Z$: 
%\[
%\sum_{t=t_1}^{t_2-1} \| x(t+1) - x (t) \| \leq c_0 + \eps ( t_2 - t_1 ) , \;
%t_1 > t_2 .
%\]

\section{The Setup}

In this paper
we start with an $n^{th}$ order linear time-invariant discrete-time plant
given by
\begin{eqnarray}
y(t+1) &=& - \sum_{i=0}^{n-1} a_{i+1} y(t-i) + \nonumber \\
&& \sum_{i=0}^{n-1} b_{i+1} u(t-i) + d(t)
\nonumber
\end{eqnarray}
\begin{eqnarray}
&=&
\underbrace{
\left[ \begin{array}{c}
y(t) \\ \vdots \\ y(t-n+1) \\ u(t) \\ \vdots \\ u(t-n+1)
\end{array} \right]^T}_{=: \phi (t)^T}
\underbrace{
\left[ \begin{array}{c}
-a_1 \\ \vdots \\ -a_n \\ b_1 \\ \vdots \\ b_n
\end{array} \right]}_{=: \theta^*}  + d(t) ,  \nonumber \\
&& \;\;\;\;\;\; \phi (t_0) = \phi_0 , \; t \geq t_0 ,
\label{plant}
\end{eqnarray}
with $y(t) \in \R$ the measured output,
$u(t) \in \R$ the control input,
and $d(t) \in \R$ the disturbance (or noise) input.
We assume that $\theta^*$ is unknown but belongs to a 
known set ${\cal S} \subset \R^{2n}$.
Associated with this plant model are the polynomials
\[
A(z^{-1} ) := 1 + a_1 z^{-1} + \cdots + a_n z^{-n} , \;\;
\]
\[
B(z^{-1} ) :=  b_1 z^{-1} + \cdots + b_n z^{-n} 
\]
and the transfer function $\frac{B(z^{-1} )}{A(z^{-1} )}$.

\begin{remark}
It is straight-forward to verify that if the system has a disturbance at both the input and output, then it can be converted to a system of the above form. 
%To see this,
%suppose that we start with the model
%\begin{eqnarray}
%y(t+1) &=& - \sum_{i=0}^{n-1} a_{i+1} y(t-i) + \sum_{i=0}^{n-1} b_{i+1} u(t-i) \\
%y_n(t) &=& y(t) + d_1(t) \\
%u(t) &=& u_n(t) + d_2(t);
%\end{eqnarray}
%here $y_n(t)$ is the measured output,
%$u_n(t)$ is the input generated by the controller,
%$d_1(t)$ is an output disturbance/noise and $d_2(t)$ is an input disturbance/noise.
%This can be rewritten as
%\begin{eqnarray*}
%y_n(t+1) &=&  - \sum_{i=0}^{n-1} a_{i+1} y_n(t-i) + \sum_{i=0}^{n-1} b_{i+1} 
%u_n(t-i) + \\
%&& 
%\underbrace{d_1 (t+1) + \sum_{i=0}^{n-1} a_{i+1} d_1(t-i) +
%\sum_{i=0}^{n-1} b_{i+1} d_2(t-i)}_{=: d(t)} ,
%\end{eqnarray*}
%which is exactly of the form (\ref{plant}).
\end{remark}

We impose an assumption on the set of admissible plant parameters.

\noindent
\framebox[85mm][l]{
\parbox{82mm}{
{\bf Assumption 1:}
${\cal S} $ is convex and compact, and for each $\theta^* \in {\cal S}$,
the 
corresponding pair of polynomials 
$A(z^{-1} ) $ and $B(z^{-1} )$ are coprime.
}}

The convexity part of the above assumption is
common in a branch of the adaptive control literature - it is used to
facilitate parameter projection, e.g.  see \cite{goodwinsin}.
The boundedness part is less common, but it is quite
reasonable in practical situations; it is used here to ensure that
we can prove uniform bounds and decay rates on the closed-loop behaviour.

The main goal here is to prove a form of stability, with a secondary
goal that of asymptotic tracking
of an 
exogenous reference signal $y^* (t)$;
since the plant may be non-minimum phase, there are limits on how well
the plant can be made to track $y^*(t)$.
To proceed we use a parameter estimator together with an adaptive
pole placement
control law.
At this point, we discuss the most critical aspect - the parameter estimator.

\subsection{Parameter Estimation}

We can write the plant as
\[
y(t+1) = \phi (t)^T \theta^* + d(t) .
\]
Given an estimate $\hat{\theta} (t)$ of $\theta^*$ at time $t$,
we define
the {\bf prediction error} by
\[
e(t+1) := y(t+1) -  \phi (t)^T  \hat{\theta} (t) ;
\]
this is a measure of the error in $\hat{\theta} (t)$.
The common way to
obtain
a new estimate is
from
the solution of the optimization problem
\[
argmin_{\theta} \{ \| \theta - \hat{\theta} (t) \| : y(t+1) = \phi (t)^T {\theta} \} ,
\]
yielding the
{\bf ideal (projection) algorithm}
\be
\that (t+1) =
\left\{
\begin{array}{ll}
\that (t) & \mbox{ if $\phi (t) = 0 $} \\
\that (t) + \frac{\phi (t)}{\phi (t)^T \phi (t)}
\,  e (t+1)  & \mbox{ otherwise.}
\end{array}
\right.
\label{orig}
\ee
Of course, if $\phi (t)$ is close to zero, numerical problems can occur, so it is the norm
in the literature (e.g. \cite{goodwin1980} and \cite{goodwinsin}) to replace this by
the following {\bf classical algorithm}:
with
$0 < \alpha < 2$ and $\beta > 0$, define
\vspace{-0.1cm}
\be
\that (t+1) = \that (t) + \frac{\alpha \phi (t)}{ \beta + \phi (t)^T \phi (t)}
e(t+1) .\footnote{It is common to make this more general by letting
$\alpha$ be time-varying.}
\label{new}
\ee
This latter algorithm is widely used, and plays a role in many discrete-time
adaptive control algorithms;
however, when this algorithm is used,
all of the results are
asymptotic, and
exponential stability and a bounded gain on the noise are never proven.
It is not hard to guess why - a careful look at the estimator shows that the gain
on the update law is small if $\phi (t)$ is small.
A more mathematically detailed argument is given in the following example.

\begin{remark}
Consider the simple first order plant
\[
y(t+1)= - a_1 y(t) + b_1 u(t) + d(t)
\]
with $a_1 \in [-2, -1]$ and $b_1 \in [1,2]$. For simplicity, we assume that 
in the estimator (\ref{new}) we
have $\alpha = \beta =1$, and, as in \cite{ydstie}, \cite{ydstie2}, \cite{naik}, \cite{wenhill}, \cite{wen} and \cite{hanfu},
we use projection to keep the parameters
estimates inside ${\cal S}$ so as to
guarantee a bounded-input bounded-state property.
Further suppose
$y^*=d=0$, and
that a
classical 
pole placement adaptive controller places the closed-loop pole at zero:
$u(t) = {\frac{\hat{a}_1 (t)}{ \hat{b}_1 (t)}} y(t) =: \hat{f} (t) y(t) $.
Suppose that
\[
y(0)  = y_0 = \eps \in (0,1),
\]
\[
\hat{\theta} (0) = \tbo{- \hat{a}_1 (0)}{ \hat{b}_1 (0)} = \tbo{1}{2}, \;
\theta^* = \tbo{2}{1}
\]
so that
$\hat{f} (0) = - 0.5$ and
$-a_1 + b_1 \hat{f} (0 ) = 1.5$, i.e. the system is initially unstable.
An easy calculation verifies that
$\hat{f} (t) \in [-2, - 0.5]$ and $
- a_1 + b_1 \hat{f} (t) 
\in [ 0,  1.5  ] $ for $t \geq 0$,
which leads to a 
crude bound on the closed loop behaviour:
$| y(t) | \leq  (1.5)^t \eps $
for $t \geq 0$.
With
$N( \eps ) := \mbox{int} [ \frac{1}{2 \ln (1.5)} \ln ( \frac{1}{\eps } ) ]$,
it follows that
\[
| y(t) | \leq \eps^{1/2} , \;\; t \in [0, N( \eps ) ] .
\]
A careful examination of the parameter estimator shows that
\[
\| \hat{\theta} (t) - \theta_0 \| \leq 10 (2)^{1/2} \eps , \;\; t \in [0, N( \eps ) ] .
\]
From the form of $\hat{f} (t)$, it follows that for small $\eps$ we have
$| -a + b_1 \hat{f} (t) | \geq 1.25 $ for $ t \in [0, N( \eps ) ]$,
in which case
\[
| y (N(\eps) ) |  \geq ( 1.25 )^{N(\eps)} \eps \;\; \Rightarrow \;\;
\left| \frac{y ( N( \eps ))}{\eps } \right| \geq (1.25)^{N(\eps )} ;
\]
since $N(\eps ) \rightarrow \infty$ as $\eps \rightarrow 0$, we see that exponential
stability is unachievable.
A similar kind of analysis can be used to prove that a bounded
gain on the noise is not achievable either.
\end{remark}

Now we return to the problem as hand - analysing the ideal algorithm
(\ref{orig}).
We will be using the ideal algorithm with projection to ensure
that the estimate remains in ${\cal S}$ for all time.
With an initial condition of $\hat{\theta} (t_0) = \theta_0 \in {\cal S}$,
for $t \geq t_0$ we
set
\be
\check{\theta} (t+1) =
\left\{ \begin{array}{ll}
\hat{\theta} (t) & 
 \;\;\mbox{ if $\phi (t) = 0 $} \\
\that (t) + \frac{\phi (t)}{\phi (t)^T \phi (t)}
\,  e (t+1) &
\;\; \mbox{otherwise,}
\end{array} \right.
\label{orig2}
\ee
which we then project onto ${\cal S}$:
\be
\hat{\theta} (t+1):= \pi_{\cal S} ( \check{\theta} (t+1) ).
\label{est2}
\ee
Because of the closed and convex property of ${\cal S}$, the projection function
is well-defined; furthermore, it has the nice property that, for every
$\theta \in \R^{2n}$ and every $\theta^* \in {\cal S}$, we have
\[
\| \pi_{\cal S} ( \theta ) - \theta^* \| \leq \| \theta - \theta^* \| ,
\]
i.e. projecting $\theta$ onto ${\cal S}$ never makes it further away from
the quantity $\theta^*$.

\subsection{Revised Parameter Estimation}

Some readers may be concerned that the original problem of dividing by a number
close to zero, which motivates
the use of classical algorithm, remains. 
Of course, this is balanced against the soon-to-be-proved benefit of
using (\ref{orig2})-(\ref{est2}).
We propose a middle ground as follows. A straight-forward analysis of $e(t+1)$ reveals that
\[
e(t+1) = - \phi (t)^T [ \hat{\theta} (t) - \theta^* ] + d (t) ,
\]
which means that 
\[
| e(t+1) | \leq 2 \| {\cal S} \| \times \| \phi (t) \| + | d(t) | .
\]
Therefore, if
\[
| e(t+1) | > 2 \| {\cal S} \| \times \| \phi (t) \| ,
\]
then the update to $\hat{\theta} (t)$ will be greater than $2 \| {\cal S} \|$,
which means that there is little information content in $e(t+1)$ - it is dominated by the disturbance. With this as motivation, and with
$\delta \in (0, \infty ]$, let us replace
(\ref{orig2}) with
\be
\check{\theta} (t+1) =
\left\{ \begin{array}{l}
\that (t) + \frac{\phi (t)}{\phi (t)^T \phi (t)}
\,  e (t+1) \\
 \;\;\;\;\;\mbox{ if $| e(t+1) | < ( 2  \| {\cal S} \| + \delta )
\| \phi (t) \|$} \\
\hat{\theta} (t) \\ \;\;\; \;\; \mbox{otherwise;}
\end{array} \right.
\label{orig3}
\ee
in the case of $\delta = \infty$, we will adopt the understanding that
$ \infty \times 0  = 0 $, in which case the above formula collapses into the
original one (\ref{orig2}).
In the case that $\delta < \infty$, we can be assured that the update term
is bounded above by $2  \| {\cal S} \| + \delta $, which should alleviate concern
about having infinite gain.
We would now
like to rewrite the update to make it more concise. To this end, we
now define 
$\rhod : \R^{2n} \times \R  \rightarrow \{ 0,1 \}$ by
\[
\rhod ( \phi (t) , e (t+1)  ) := 
\]
\[
\left\{ \begin{array}{ll}
1 & \;\; \mbox{ if } | e(t+1) | < ( 2  \| {\cal S} \| + \delta ) 
\| \phi (t) \| \\
0 & \;\; \mbox{otherwise, } 
\end{array} \right.
\]
yielding a more concise way to write the estimation algorithm update:
\be
\check{\theta} (t+1) = \hat{\theta} (t) + \rhot 
\frac{ \phi (t)}{ \phi (t)^T  \phi (t)} e(t+1) ;
\label{estimator1}
\ee
once again, we project this onto ${\cal S}$:
\be
\hat{\theta} (t+1):= \pi_{\cal S} ( \check{\theta} (t+1) ).
\label{estimator2}
\ee

\subsection{Properties of the Estimation Algorithm}

Analysing the closed-loop system will require a careful
analysis of the estimation algorithm.
We define
the
parameter estimation error by
\[
\tilde{\theta} (t): = \hat{\theta} (t) - \theta^* ,
\]
and the corresponding Lyapunov function associated with
$\tilde{\theta} (t)$, namely
$V(t) : = \tilde{\theta} (t)^T \tilde{\theta} (t) $.
In the following result we list a property of $V(t)$; it is a generalization
of what is well-known for the classical algorithm (\ref{new}).
%
%\noindent
%\framebox[160mm][l]{
%\parbox{155mm}{
\noindent
\framebox[85mm][l]{
\parbox{82mm}{
\begin{prop}
For every $t_0 \in \Z$, $\phi_0 \in \R^{2n}$, ${\theta}_0 \in {\cal S}$,
$\theta^* \in {\cal S}$,
$d \in \linf$, and $\delta \in ( 0, \infty ]$, when the estimator 
(\ref{estimator1}) and (\ref{estimator2})
is applied to the plant (\ref{plant}),
the following holds:
\be
 \| \hat{\theta} (t+1) - \hat{\theta}  (t) \|
\leq   \rhot \frac{ |e(t+1)| }{ \|  \phi (t) \| } ,
\; t \geq t_0 ,
\label{errorsum}
\ee
\begin{eqnarray*}
V(t)
& \leq &
V( t_0 ) +   \sum_{j=t_0}^{t-1}
\rhoj \times
\end{eqnarray*}
\[
\;\;\;
[ -\frac{1}{2} 
\frac{[e (j+1) ]^2}{ \| \phi (j) \|^2} + 
2 \frac{[ d(j)]^2}{ \| \phi (j) \|^2}] , \;\; t \geq t_0+1 .
\]
\end{prop}
}}

\vspace{0.3cm}
\noindent
{\bf Proof:}
See the Appendix.
$\Box$

\subsection{The Control Law}
The elements of $\hat{\theta} (t)$ are partitioned in a natural way as
\[
\hat{\theta} (t) =
\left[
\begin{array}{cccccc}
- \hat{a}_1 (t)  &
\cdots &
- \hat{a}_n (t) &
 \hat{b}_1 (t) &
\cdots &
 \hat{b}_n (t) 
\end{array}
\right]^T
\]
Associated with $\hat{\theta} (t) $ are the polynomials
\[
\hat{A} (t, z^{-1} ) := 1 + \hat{a}_1 (t) z^{-1} + \cdots + 
\hat{a}_n (t) z^{-n} , \;\;
\]
\[
\hat{B}(t,z^{-1} ) :=  \hat{b}_1 (t) z^{-1} + \cdots + \hat{b}_n (t) z^{-n} .
\]
While we can use an $n-1^{th}$ order {\bf proper} controller to carry out pole placement, it will be convenient to follow the lead of 
\cite{wenhill} and use an $n^{th}$ order {\bf strictly
proper} controller.
In particular, we first choose a $2n^{th}$ order monic polynomial
\[
A^* (z^{-1} ) = 1 + a_1^* z^{-1} + \cdots + a_{2n}^* z^{-2n} 
\]
so that $z^{2n} A^* (z^{-1} )$ has all of its zeros in $\D^o$.
Next, we choose two polynomial 
\[
\hat{L} (t, z^{-1} ) = 1 + \hat{l}_1 (t) z^{-1} + \cdots + \hat{l}_n (t)
z^{-n} 
\]
and
\[
\hat{P} (t,z^{-1} ) = \hat{p}_1 (t) z^{-1} + \cdots + \hat{p}_n (t) z^{-n} 
\]
which satisfy the equation
\be
\hat{A} ( t , z^{-1} ) \hat{L} (t, z^{-1} ) +
\hat{B} ( t , z^{-1} ) \hat{P} (t,z^{-1} ) = A^* ( z^{-1} ) ;
\label{poleplace}
\ee
given the assumption that the $\hat{A} ( t , z^{-1} )$ and
$\hat{B} ( t , z^{-1} )$ are coprime, it is well known that
there exist {\bf unique}
$\hat{L} (t, z^{-1} )$ and $\hat{P} (t,z^{-1} )$ which satisfy this
equation.
Indeed, it is easy to prove that the coefficients of 
$\hat{L} (t, z^{-1} )$ and $\hat{P} (t,z^{-1} )$ are analytic
functions of $\hat{\theta} (t) \in {\cal S}$.
%Since ${\cal S}$ is compact, this means that
%\[
%\sup_{\hat{\theta} (t)} \sum_{i=1}^n ( | \hat{l}_i (t) | + | \hat{p}_i (t) | )
%< \infty .
%\]

In our setup we have an exogenous signal $y^* (t)$.
At time $t$ we choose $u(t)$ so that
\[
 u(t) =  - \hat{l}_1  (t-1) u(t-1) - \cdots - 
\]
\[
- \hat{l}_n (t-1) u(t-n) 
 - \hat{p}_1 (t-1) [ y(t-1)- y^* (t-1) ] - \cdots 
\]
\be
- \hat{p}_n (t-1) [ y(t-n)- y^* (t-n) ]  .
\label{conguy}
\ee
So the overall controller consists of the estimator (\ref{estimator1})-(\ref{estimator2})
together with (\ref{conguy}).\footnote{We also implicitly
use a pole placement procedure to
obtain the controller parameters from the plant parameter estimates;
this entails solving a linear equation.}

It turns out that 
we can write down a state-space model of our closed-loop system with
$\phi (t) \in \R^{2n}$ as the state. Only two elements of $\phi$ have a complicated
description:
\begin{eqnarray*}
\phi_1 (t+1) &=& y(t+1) 
= e(t+1) + \hat{\theta} (t)^T \phi (t) ,
\end{eqnarray*}
%\begin{eqnarray*}
%\phi_2 (t+1) &=& y(t) = \phi_1 (t) \\
%& \vdots&  \\
%\phi_n (t+1) &=& y(t-n+2) = \phi_{n-1} (t)
%%\end{eqnarray*}
\begin{eqnarray*}
\phi_{n+1}(t+1) &=& u(t+1) \\
&=&
- \sum_{i=1}^{n} \{ \hat{l}_i  (t) u(t+1-i) - \\
&&
 \hat{p}_i (t) [ y(t+1-i) - y^* (t+1-i) ] \}  
\end{eqnarray*}
\[
=
\left[ \begin{array}{cccccc}
- \hat{l}_1 (t) & \cdots & - \hat{l}_n (t) & - \hat{p}_1 (t) & \cdots & - \hat{p}_n (t)
\end{array} \right] \phi (t) +  
\]
\[
 \;\;\;
\sum_{i=1}^n \hat{p}_i (t) y^* ( t+1 - i )  .
\]
With $e_i \in \R^{2n}$ the $i^{th}$ normal vector,
if we now define
\[
\bar{A} (t) :=
\left[ \begin{array}{cccc}
-\hat{a}_1 (t) & -\hat{a}_2 (t) & \cdots & -\hat{a}_n (t) \\
1 & 0 & \cdots & 0  \\
 &  \ddots & &  \vdots  \\
 & & 1 & 0  \\
- \hat{p}_1 (t) & - \hat{p}_2 (t) & \cdots & -\hat{p}_n (t) \\
0 & \cdots & \cdots & 0 \\
\vdots & \cdots & \cdots & \vdots \\
0 & \cdots & \cdots & 0  \\
\end{array} \right. 
\]
\[
\left. \begin{array}{cccc}
\hat{b}_1 (t) & \cdots & \cdots & \hat{b}_n (t) \\
0 & \cdots & \cdots & 0 \\
\vdots  & \cdots & \cdots & \vdots \\
 0  & \cdots & \cdots & 0 \\
- \hat{l}_1 (t) & -\hat{l}_2 (t) & \cdots & - \hat{l}_n (t) \\
1 & 0 & \cdots & 0 \\
 & \ddots &  & \vdots \\
 &  & 1 & 0 \\
\end{array} \right] ,
\]
\be
B_1:= 
e_1 , 
%\left[ \begin{array}{c}
%1 \\ 0 \\ \vdots \\ 0 \end{array} \right] ,
\;\;
B_2 := e_{n+1} ,
%  \left[ \begin{array}{c}
%0 \\ \vdots \\ 0 \\ 1 \\ 0 \\ \vdots \\ 0
%\end{array} \right] ,
\;\;
r(t) := \sum_{i=1}^{n} \hat{p}_i (t) y^* (t+1-i)  ,
\label{rdef}
\ee
then 
the following key equation holds:
\be
\phi (t+1) = \bar{A} (t) \phi (t) + B_1  e(t+1)  + B_2 r(t) ;
\label{keyeq}
\ee
notice that the characteristic equation of $\bar{A} (t)$
always equals
$z^{2n} A^* ( z^{-1} )$.
%(check ***$r(t)$ definition ***)
%Using (\ref{predict}), this can be rewritten as
%\be
%\phi (t+1) = [ \bar{A} (t) - B_1 \tilde{\theta} (t)^T ] \phi (t) + B_1 d(t) +   B_2 r(t) .
%\label{neweq}
%\ee
Before proceeding, define
\[
\bar{a} := \max \{ \| \bar{A} ( \hat{\theta} ) \| : \hat{\theta} \in {\cal S} \} .
\]

\section{Preliminary Analysis}

The closed-loop system given in 
(\ref{keyeq}) arises in classical adaptive control approaches in slightly
modified fashion, so we will borrow some tools from there. More specifically,
the following result was proven by Kreisselmeir \cite{kreiss2},
in the context of proving that a slowly time-varying adaptive control
system is stable (in a weak sense); we are providing a special case of his
technical lemma
to minimize complexity.\footnote{Furthermore, in
\cite{kreiss2} it is assumed that $\alpha_i $ and
$\beta_i $ are strictly greater than zero, but it is trivial to extend this
to allow for zero as well.}

\noindent
%\framebox[160mm][l]{
%\parbox{155mm}{
\noindent
\framebox[85mm][l]{
\parbox{82mm}{
\begin{prop} \cite{kreiss2}
Consider the discrete-time system
\[
x (t+1) = [ A_{nom} (t) + \Delta (t) ] x(t)
\]
with $\Phi (t, \tau )$ denoting the 
corresponding state transition matrix.
Suppose that there exist constants
$\sigma \in (0,1)$, $\gamma_1 >1$, $\alpha_i \geq 0$, and $\beta_i \geq 0$
so that

\vspace{0.2cm}
\noindent
{(i)}
for all $t \geq t_0 $, we have
$\| A_{nom} (t)^i \| \leq \gamma_1 \sigma^i , \; i \geq 0 $;

\vspace{0.2cm}
\noindent
{(ii)}
for all $t > \tau $ we have
\[
\sum_{i = \tau}^{t-1} \| A_{nom} (i+1) - A_{nom} (i) \| \leq
\]
\[
\;\;\;\;\;
\alpha_0 + \alpha_1 ( t- \tau )^{1/2} + \alpha_2 ( t- \tau ) 
\]
and
$\sum_{i = \tau}^{t-1} \|  \Delta (i) \| \leq
\beta_0 + \beta_1 ( t- \tau )^{1/2} + \beta_2 ( t- \tau ) $;

\vspace{0.2cm}
\noindent
{(iii)}
there exists a $\mu \in ( \sigma , 1 )$ and $N \in \N$ satisfying
$\alpha_2 + \frac{\beta_2}{N} < \frac{1}{N \gamma_1} ( \frac{\mu}{\gamma_1^{1/N}} -
\sigma ) $.

\vspace{0.2cm}
\noindent
Then
there exists a constant $\gamma_2$ so that 
the transition matrix satisfies
\[
\| \Phi (t, \tau  ) \| \leq \gamma_2 \mu^{t - \tau }  , \; t \geq \tau .
\]
\end{prop}
}}

\begin{remark}
We apply the above proposition in the following way.
Suppose that $\sigma \in (0,1)$, $\gamma_1 >1$, $\alpha_i \geq 0$, $\beta_i \geq 0$
are such conditions (i) and (ii) hold.
If $\mu \in ( \sigma , 1 )$, then it follows that
$\frac{\mu}{\gamma_1^{1/N}} -
\sigma  > 0 
$
for large enough $N \in \N$, so condition (iii) will hold as well as long as
$\alpha_2 $ and $\beta_2$ are small enough.
\end{remark}

In applying Proposition 2, the matrix $\bar{A} (t)$ will play the role of
$A_{nom} (t)$.
A key requirement is that Condition (i) holds: the following 
provides relevant bounds.
Before proceeding, let
\[
\underline{\lambda} \; := \; \max \{ | \lambda | : \lambda  \;  \mbox{is a root of $z^{2n} A^* (z^{-1})$} \}.
\]

\noindent
%\framebox[160mm][l]{
%\parbox{155mm}{
\noindent
\framebox[85mm][l]{
\parbox{82mm}{
\begin{lemma}
For every $\delta \in ( 0 , \infty ]$ and $\sigma \in (\underline{\lambda} ,1 )$
there exists a constant 
$\gamma \geq 1$ so that for every $t_0 \in \Z$,
$\theta_0 \in {\cal S}$, $\hat{\theta}^* \in {\cal S}$,  and $y^*,d \in \linf$, 
when the controller 
(\ref{estimator1}), (\ref{estimator2})
and
(\ref{conguy}) is applied to the plant (\ref{plant}), 
the matrix $\bar{A} (t)$ satisfies, for every $t \geq t_0$:
\[
\| \bar{A} (t)^k \| \leq \gamma \sigma^k , \; k \geq 0,
\]
and
for every $t > k \geq t_0$:
\[
\sum_{j=k}^{t-1} \| \bar{A} (j+1) - \bar{A} (j) \| \leq \gamma \times
\]
\[
[\sum_{j=k}^{t-1} \rhoj  \frac{e(j+1)^2}{\| \phi (j)\|^2} ]^{1/2} ( t-k)^{1/2} .
\]
\end{lemma}
}}

\noindent
{\bf Proof:}  See the Appendix.
$\Box$

\section{The Main Result}

%Before proceeding, recall
%that $\underline{\lambda} := $ the magnitude of the largest root of $z^{2n} A^* (z^{-1})$.

\noindent
\framebox[85mm][l]{
\parbox{82mm}{
\begin{theorem}
For every
$\delta \in ( 0 , \infty ]$ and
$\lambda \in ( \underline{\lambda} ,1)$ there exists a $c>0$ so that for
every $t_0 \in \Z$, $\theta_0 \in {\cal S}$, ${\theta}^*  \in {\cal S}$,
$\phi_0 \in \R^{2n}$,
and $y^*, d \in \ell_{\infty}$,  when the adaptive controller 
(\ref{estimator1}), (\ref{estimator2})
and
(\ref{conguy}) is applied
to the plant (\ref{plant}), the following bound holds:
\[
\| \phi (k) \| \leq c \lambda^{k- t_0} \| \phi_0 \| +
\]
\be
\sum_{j=t_0}^{k-1}  c \lambda^{k-1-j} {( | r (j) | + | d(j)| )} , \;\;
k \geq  t_0 .
\label{thm1_eq}
\ee
\end{theorem}
}}

\begin{remark}
We see from 
(\ref{rdef}) that $r(t)$ is a weighted sum of
$\{ y^* (t) ,..., y^* ( t-n+1 ) \}$. Hence, there exists a
constant $\bar{c}$ so that the bound (\ref{thm1_eq})
can be rewritten as
\[
\| \phi (k) \| \leq c \lambda^{k- t_0} \| \phi_0 \| +
\sum_{j=t_0}^{k-1}  c \lambda^{k-1-j}  | d (j) | + 
\]
\[
\;\;\;
\sum_{j=t_0-n+1}^{k-1}  \bar{c} \lambda^{k-1-j} | y^*(j)|  , \;\;
k \geq  t_0 .
\]
\end{remark}

\begin{remark}
Theorem 1 implies that the system has a bounded gain (from 
$d$ and $r$ to $y$) in every $p-$norm.
More specifically, for $p= \infty$ we see immediately from (\ref{thm1_eq})
that
\[
\| \phi (k) \| \leq c  \| \phi_0 \| +
\frac{c}{1- \lambda} \sup_{ \tau \in [t_0 , k] } [ | r( \tau ) | +
| d( \tau ) | ) ] , \;\; k \geq t_0 .
\]
Furthermore, for $1 \leq p < \infty$ it follows from
Young's Inequality applied to (\ref{thm1_eq}) that
\[
[ \sum_{j=t_0}^k \| \phi (j) \|^p ]^{1/p} \leq
\frac{c}{(1 - \lambda^p)^{1/p}} \| \phi_0 \| +
\]
\[
\frac{c}{1 - \lambda} \{
[ \sum_{j=t_0}^k \| r (j) \|^p ]^{1/p} +
[ \sum_{j=t_0}^k \| d (j) \|^p ]^{1/p} \} , \; k \geq t_0 .
\]
\end{remark}

\begin{remark}
Most pole placement adaptive controllers are proven to yield a weak form of stability,
such as boundedness (in the presence of a non-zero disturbance) or asymptotic
stability (in the case of a zero disturbance), which means that details surrounding initial conditions
can be ignored.
Here the goal is to prove a stronger, linear-like, convolution
bound, so it requires more detailed analysis.
\end{remark}

\begin{remark}
With $\hat{G} (t, z^{-1} ) =\sum_{i=1}^{2n}
\hat{g}_i (t) z^{-i} := \hat{B} (t, z^{-1} ) \hat{P} (t,z^{-1} )$
it is possible to use arguments like those in \cite{goodwinsin} to prove,
when the disturbance $d$ is identically zero, a weak tracking result of the form
\[
\lim_{t \rightarrow \infty} [
\sum_{i=0}^{2n} a_i^* y(t-i) -
\sum_{i=1}^{2n} \hat{g}_i (t) y^* ( t-i) ] = 0 .
\]
Since the main goal of the paper is on stability issues, we omit the
proof. However, we do discuss step tracking in a later section.
\end{remark}

\vspace{0.5cm}
\noindent
{\bf Proof:}
Fix $\delta \in (0, \infty ]$ and $\lambda \in ( \underline{\lambda}, 1)$. Let $t_0 \in \Z$, 
$\theta_0 \in {\cal S}$,
$\theta^* \in {\cal S}$,
$\phi_0 \in \R^{2n}$,
and 
$y^*, d\in \linf$ be arbitrary.
Define $r$ via (\ref{rdef}).
Now choose
$\lambda_1 \in 
( \underline{\lambda} ,  \lambda )$.

We have to be careful in how to apply Proposition 2 to (\ref{keyeq}) -
we need the $\Delta (t)$ term to be something which
we can bound using Proposition 1.
So
define
\be
\Delta (t):= \rhot \frac{e(t+1)}{\| \phi (t) \|^2} B_1 \phi (t)^T ;
\label{ddef}
\ee
it is easy to check that
\[
\Delta (t) \phi (t) = \rhot B_1 e(t+1)
\]
and that
\[
\| \Delta (t) \| = \rhot \frac{| e(t+1)|}{\| \phi (t) \| } ,
\]
which is a term which plays a key role in Proposition 1. We can now rewrite
(\ref{keyeq}) as
\[
\phi (t+1) = [ \bar{A} (t) + \Delta (t) ] \phi (t) +
\]
\be
 B_1  \underbrace{[ 1 - 
\rhot ]
 e(t+1) 
 }_{=: \eta (t)}  +   B_2 r(t) .
\label{neweq}
\ee
If $\rhot =1$ then $\eta (t) =0$, but if $\rhot  = 0$ then
\[
| e(t+1) | \geq ( 2 \|  {\cal S} \| + \delta ) \| \phi (t) \| ;
\]
but
we also know that
\[
e(t+1) = - \tilde{\theta} (t) \phi (t) + d(t) 
\]
\be
\; \Rightarrow \; | e(t+1) | \leq 2 \| {\cal S} \| \times \| \phi (t) \| + | d(t) | ;
\label{ebd}
\ee
combining these equations we have
\[
 ( 2 \|  {\cal S} \| + \delta ) \| \phi (t) \|  \leq 2 \| {\cal S} \| \times \| \phi (t) \| + | d(t) | ,
\]
which implies that
$\| \phi (t) \| \leq \frac{1}{\delta} | d(t) | $;
it is easy to check that this holds even when $\delta = \infty$. Using (\ref{ebd}) we conclude that
\be
| \eta (t) | \leq ( \frac{2 \| {\cal S} \|}{\delta} +1 ) | d(t) | , \; t \geq t_0 .
\label{etabd}
\ee

We now analyse (\ref{neweq}). 
We let $\Phi (t, \tau )$ denote the transition matrix associated with
$ \bar{A} (t) + \Delta (t)$; this matrix clearly implicitly depends on
$\theta_0$, $\theta^*$, $d$ and $r$.
From Lemma 1 there exists a constant $\gamma_1$ so that
\be
\| \bar{A} (t)^i \| \leq \gamma_1 \lambda_1^i, \;\; i \geq 0, \; t \geq t_0 ,
\label{babd}
\ee
and for every $t > k \geq t_0$, we have
\[
\sum_{j=k}^{t-1} \| \bar{A} (j+1) - \bar{A} (j) \|  \leq 
\]
\be
  \gamma_1 [  \sum_{j=k}^{t-1}  \rhojm
\frac{ | e(j+1)|^2}{ \| \phi (j) \|^2 } ] ^{1/2} ( t-k)^{1/2}  .
\label{babdd}
\ee
Using the definition of $\Delta$ given in (\ref{ddef}) and the Cauchy-Schwarz inequality we also have
\[
\sum_{j=k}^{t-1} \|  \Delta  (j) \|  \leq  
 [  \sum_{j=k}^{t-1}  \rhojm \times
\]
\be
\frac{ | e(j+1)|^2}{ \| \phi (j) \|^2 } ] ^{1/2} ( t-k)^{1/2} , \;
t > k \geq t_0 .
\label{dbdd}
\ee
At this point we consider two cases: the easier case in which there is no noise,
and the harder case in which there is noise.

\noindent
{\bf Case 1}: $d(t) = 0$, $t \geq t_0$.

Using the bound on $\eta (t)$ given in (\ref{etabd}),
in this case (\ref{neweq}) becomes
\be
\phi (t+1) = [ \bar{A} (t)  + 
\Delta (t)]  \phi (t) +  B_2 r(t), \; t \geq t_0 .
\label{deq}
\ee
The bound on $V(t)$ given by Proposition 1
simplifies to
\[
V(t) \leq V ( t_0  ) - \frac{1}{2} \sum_{j=t_0}^{t-1} \rhojm  \frac{ [e(j+1)]^2}
{ \| \phi (j) \|^2} ,
\]
\[
\;\;\;\; \; t \geq t_0 +1 .
\]
Since $V( \cdot ) \geq 0$ and $ V( t_0  ) = \| \theta_0 - \theta^* \|^2
\leq 4 \| {\cal S} \|^2$, this means that
\[
\sum_{j=t_0}^{t-1} \rhojm   \frac{ [e(j+1)]^2}
{\|  \phi (j)\|^2} \leq  2 V( t_0  )  \leq 8 \| {\cal S} \| ^2.
\]
Hence, 
from (\ref{babdd}) and (\ref{dbdd}) we conclude that
\[
\sum_{j=k}^{t-1} \| \bar{A} (j+1) - \bar{A} (j) \|  \leq
% \gamma_1 [  \sum_{j=k}^{t-1}  \rhojm
%\frac{ | e(j+1)|^2}{ \| \phi (j) \|^2 } ] ^{1/2} ( t-k)^{1/2} \leq 
8^{1/2}  \gamma_1 \| {\cal S} \|
( t-k)^{1/2} ,
\]
\[
\sum_{j=k}^{t-1} \| \Delta (j)  \|  \leq
% [  \sum_{j=k}^{t-1}  \rhojm
%\frac{ | e(j+1)|^2}{ \| \phi (j) \|^2 } ] ^{1/2} ( t-k)^{1/2} \leq 
8^{1/2} \| {\cal S} \| 
( t-k)^{1/2} 
 , \; t > k \geq t_0 .
\label{dbdd2}
\]
Now we apply Proposition 2: we set
\[
\alpha_0 = \beta_0 = \alpha_2 = \beta_2 = 0 ,
\]
\[
 \alpha_1 =   8^{1/2} \gamma_1 \| {\cal S} \| ,
\; \beta_1 = 8^{1/2} \| {\cal S} \| ,
\; \mu = \lambda.
\]
Following Remark 3 it
is now trivial to choose $N \in \N$ so that
$\frac{\lambda}{ \gamma_1^{1/N}} - \lambda_1 > 0 $,
namely
\be
N = \mbox{int} [ \frac{\ln ( \gamma_1 )}{\ln ( \lambda ) - \ln ( \lambda_1 ) } ] + 1 ,
\label{ndef}
\ee
which means that
\[
0 = \alpha_2 + \frac{\beta_2}{N} < \frac{1}{N \gamma_1 } (  \frac{\lambda}{ \gamma_1^{1/N}} - \lambda_1 ) .
\]
From Proposition 2 we see that there exists a constant $\gamma_2$ so that the state
transition matrix $\Phi ( t, \tau )$ corresponding to $\bar{A} (t) + \Delta (t)$ satisfies
\[
\| \Phi (t, \tau ) \| \leq \gamma_2 \lambda^{t - \tau } , \;
t \geq \tau \geq  t_0 .
\]
If we now apply this to
(\ref{deq}), we end up with the desired bound:
\[
\| \phi (k) \| \leq \gamma_2 \lambda^{k- t_0} \| \phi ( t_0) \| +
\sum_{j=t_0}^{k-1}  \gamma_2 \lambda^{k-1-j} | r (j) |  , \;\;
k \geq  t_0 .
\]

\noindent
{\bf Case 2:} $d(t) \neq 0$ for some $t \geq t_0$.

This case is much more involved since noise can radically
affect parameter estimation. Indeed, even if the parameter estimate is
quite accurate at a point in time, the introduction of a large noise signal
(large relative to the size of $\phi (t)$) can create a highly inaccurate
parameter estimate.
To proceed we partition the timeline into two parts: one in which the noise is
small versus $\phi$ and one where it is not; the actual choice of the
line of division will become clear
as the proof progresses.
To this end, with
$\eps >0 $ to be chosen shortly,
partition $\{ j \in \Z : j \geq t_0 \}$ into two sets:
\[
S_{good} :=  \{ j \geq t_0 : \phi (j) \neq 0 \mbox{ and }
 \frac{[ d(j)]^2}{ \| \phi (j) \|^2}  < \eps \} ,
\]
\[
S_{bad} := \{ j \geq t_0 : \phi (j) = 0 \mbox{ or }
 \frac{[ d(j)]^2}{ \| \phi (j) \|^2} \geq \eps \} ;
\]
clearly
$\{ j \in \Z : \;\; j \geq t_0 \} = S_{good} \cup S_{bad}  $.
Observe that this partition clearly depends on $\theta_0$, $\theta^*$, $\phi_0$,
$d$ and $r/y^*$.
We will apply Proposition 2 to analyse the closed-loop system
behaviour on $S_{good}$;
on the other hand, we will easily obtain bounds on the system behaviour on
$S_{bad}$.
Before doing so, we partition the time index
$\{ j \in \Z: j \geq t_0 \}$ into
intervals which oscillate between $S_{good}$ and $S_{bad}$.
To this end, it is easy to see that we can define
a (possibly infinite) sequence of intervals of the form
$[ k_i , k_{i+1} )$ satisfing:

\vspace{0.2cm}
\noindent
{(i)}
$k_1 = t_0$, and

\vspace{0.2cm}
\noindent
{(ii)}
$[ k_i , k_{i+1} )$ either belongs to $S_{good}$ or $S_{bad}$, and

\vspace{0.2cm}
\noindent
{(iii)}
if $k_{i+1} \neq \infty$ and $[ k_i , k_{i+1} )$ belongs to $S_{good}$
(respectively, $S_{bad}$), then the interval
$[ k_{i+1} , k_{i+2} )$ must belong to $S_{bad}$ (respectively,
$S_{good}$).

Now we turn to analysing the behaviour during each interval.

\noindent
{\bf Sub-Case 2.1:}
$[ k_i , k_{i+1} )$ lies in $S_{bad}$.

Let $j \in [ k_i , k_{i+1} )$ be arbitrary.
In this case
either
$\phi (j) = 0$
or
$\frac{[ d(j)]^2}{ \| \phi (j) \|^2} \geq \eps$
holds.
In either case we have
\be
\| \phi (j) \| \leq \frac{1}{\eps^{1/2}}  | d(j) |, \; j \in [ k_i , k_{i+1} ) .
\label{phibd1}
\ee
From (\ref{keyeq}) and (\ref{ebd}) we see that 
\begin{eqnarray}
\| \phi (j+1) \| & \leq & 
\bar{a} \| \phi (j) \| + \nonumber \\
&&  ( 2 \| {\cal S} \| \times \| \phi (j) \| + | d(j) | + | r(j) | ) \nonumber
\\
%& \leq &
%( \bar{a} + 2 \| {\cal S} \|  ) \| \phi (j) \| + | d(j) | + | r(j) | \nonumber \\
& \leq & 
[ 1+ \underbrace{( \bar{a} + 2 \| {\cal S} \| )}_{=: \gamma_3} 
\frac{1}{\eps^{1/2}} ] | d(j) | + \nonumber \\
&& \;\;\; | r(j) | , \; j \in [ k_i , k_{i+1} ) .
\label{bad_bound}
\end{eqnarray}
If we combine this with (\ref{phibd1}) we conclude that
\be
\| \phi (j) \| \leq
\left\{ \begin{array}{l}
 \frac{1}{\eps^{1/2}} | d(j) |  \\ \;\;\;\;\;\;\;\; \mbox{if } j=k_i \\
( 1 + \frac{\gamma_3}{ \eps^{1/2} } )  | d(j-1) | + | r(j-1) | \\
\;\; \;\;\;\;\;\; \mbox{if } j = k_i+1,..., k_{i+1} .
\end{array}
\right.
\label{phibd}
\ee

\vspace{0.2cm}
\noindent
{\bf Sub-Case 2.2:}
$[ k_i , k_{i+1} )$ lies in $S_{good}$.

Let $j \in [ k_i , k_{i+1} )$ be arbitrary.
In this case
$\phi (j) \neq 0$
and
\[
 \frac{[ d(j)]^2}{ \| \phi (j) \|^2}  < \eps, \;\;
j \in  [k_i , k_{i+1} ) ,
\]
%This means that
%$|d (j) | < \eps^{1/2} \| \phi (j) \| $,
%so
%\[
%| e(j+1) | \leq ( 2 \| {\cal S} \| + \eps^{1/2}) \| \phi (j) \| ;
%\]
%since we have chosen $\eps^{1/2} < \delta$, this means that
%$\rhoj =1$, so we conclude that
which implies that
\be
\rhoj \frac{ d(j)^2}{\| \phi (j) \|^2} < \eps , \; j \in  [k_i , k_{i+1} ).
\label{dbd}
\ee
From Proposition 1 we have that
\[
V( \bk ) \leq V( \uk ) + \sum_{j=\uk}^{\bk-1} \rhoj  \times
\]
\[
\frac{ - \frac{1}{2} 
e(j+1)^2 + 2 d(j)^2}
{\| \phi (j) \|^2} , \; k_i \leq \uk < \bk \leq k_{i+1} ;
\]
using (\ref{dbd}) and
the fact that $0 \leq V( \cdot ) \leq 4 \| {\cal S} \|^2$, we obtain
\[
\sum_{j=\uk}^{\bk-1} \rhoj  \frac{  e(j+1)^2 }{ \| \phi (j)\|^2}  
\]
\begin{eqnarray*}
& \leq &
2 V( \uk ) + 2 \sum_{j=\uk}^{\bk-1} \rhojm  \frac{  2 d(j)^2}
{\| \phi (j)\|^2} \\
& \leq & 8 \|{\cal S } \|^2  + 4 \eps ( \bk - \uk ) , \;
k_i \leq \uk < \bk \leq k_{i+1} .
\end{eqnarray*}
Hence,
using this in
(\ref{babdd}) and (\ref{dbdd}) yields
\[
\sum_{j=\uk}^{\bk-1} \| \bar{A} (j+1) - \bar{A} (j) \|
\]
\begin{eqnarray*}
& \leq &
\gamma_1 [ 8 \|{\cal S } \|^2  + 4 \eps ( \bk - \uk )]^{1/2} ( \bk - \uk )^{1/2} \\
& \leq & \gamma_1  8^{1/2} \|{\cal S } \| (  \bk - \uk )^{1/2} + 2 \gamma_1 \eps^{1/2} 
 ( \bk - \uk  ) ,\\
&&
 \;\;\;
k_i \leq \uk  < \bk \leq k_{i+1} ,
\end{eqnarray*}
as well as
\begin{eqnarray*}
\sum_{j=\uk}^{\bk-1} \| \Delta (j) \|  & \leq &
[ 8 \|{\cal S } \|^2  + 4 \eps ( \bk - \uk )]^{1/2} ( \bk - \uk )^{1/2} \\
& \leq &
  8^{1/2} \|{\cal S } \| ( \bk - \uk) ^{1/2} + 2 \eps^{1/2}  ( \bk - \uk ) , \\
 && \;\;
k_i \leq \uk <  \bk \leq k_{i+1} .
\end{eqnarray*}
Now we will apply Proposition 2: we set
\[
\alpha_0 = \beta_0 =  0 ,
\; \alpha_1 =  \gamma_1 8^{1/2} \|{\cal S } \|   ,
\; \beta_1 =   8^{1/2} \|{\cal S } \|  ,
\;
\]
\[
\alpha_2 = 2 \gamma_1 \eps^{1/2} , \;
\beta_2 = 2 \eps^{1/2} , \;
\mu = \lambda.
\]
With 
$N$ chosen as in Case 1 via (\ref{ndef}), we have that
$\underline{\delta} := \frac{\lambda}{ \gamma_1^{1/N}} - \lambda_1 > 0 $;
we need
\[
\alpha_2 + \frac{\beta_2}{N} < \frac{1}{N \gamma_1} \underline{\delta} ,
\]
%or equivalently
%\begin{eqnarray*}
%2 \gamma_1 \eps^{1/2} + \frac{2 \eps^{1/2}}{N} & < & \frac{1}{N \gamma_1} \underline{\delta} \\
%\Leftrightarrow \;  \eps < \frac{ \underline{\delta}^2}{4 \gamma_1^2 ( \gamma_1 N + 1 )^2} ,
%\end{eqnarray*}
which will certainly be the case if we set
$\eps :=   \frac{ \underline{\delta}^2}{8 \gamma_1^2 ( \gamma_1 N + 1 )^2}  $.
From Proposition 2 we see that there exists a constant $\gamma_4$ so that the state
transition matrix $\Phi (t , \tau )$ corresponding to $\bar{A} (t) +
\Delta (t)$ satisfies
\[
\| \Phi (t, \tau ) \| \leq \gamma_4 \lambda^{t - \tau } , \;
k_i \leq \tau \leq t \leq k_{i+1} .
\]
If we now apply this to
(\ref{neweq}) and use (\ref{etabd}) to provide a bound on
$\eta (t)$, we end up with
\[
\| \phi (k) \| \leq \gamma_4 \lambda^{k- k_i} \| \phi ( k_i) \| +
( 2 \frac{\| {\cal S} \|}{\delta} +1 ) \times
\]
\be
\sum_{j=k_i}^{k-1}  \gamma_4 \lambda^{k-1-j} ( | r (j) | + | d(j) | )  , \;\;
k_{i} \leq k \leq k_{i+1} .
\label{good}
\ee
This completes Sub-Case 2.2.

Now we combine Sub-Case 2.1 and Sub-Case 2.2 into a general
bound on $\phi(t)$.
Define 
\[
\gamma_5 := \max \{ 1, 1 + \frac{\gamma_3}{\eps^{1/2}},  \gamma_4 , 
\gamma_4 (  2 \frac{\| {\cal S} \|}{\delta} +2 + \frac{\gamma_3}{\eps^{1/2}} ) \} .
\]
It remains to prove

\noindent
{\bf Claim:} The following bound holds:
\[
\| \phi (k) \| \leq \gamma_5 \lambda^{k-t_0}
\| \phi _0  \| + \sum_{j=t_0}^{k-1}
\gamma_5 \lambda^{k-1-j} ( |r(j) | + | d(j) | ) , \;\;
\]
\be 
\;\;\;\;\;\;
k \geq t_0  .
\label{nice}
\ee

\noindent
{\bf Proof of the Claim:}

If $[k_1 , k_2 ) = [t_0, k_2) \subset S_{good}$, then (\ref{nice}) holds for $k \in [t_0 , k_2 ]$ by
(\ref{good}).
If $[t_0 , k_2 ) \subset S_{bad}$, then from
(\ref{phibd}) we obtain
\[
\| \phi (j) \| \leq
\left\{ \begin{array}{l}
\| \phi ( k_1 ) \| = \| \phi_0 \|  \\
   \;\;\;\;\; j=k_1 = t_0 \\
( 1 + \frac{\gamma_3}{\eps^{1/2}})  | d(j-1) | + | r(j-1) | \\
 \;\;\; \;\; j = k_1+1,..., k_{2} ,
\end{array}
\right.
\]
which means that (\ref{nice}) holds for $k \in [t_0 , k_2 ]$ for this case as well.

We now use induction - suppose that (\ref{nice}) holds for $k \in [k_1 , k_i ]$; we need to 
prove that it holds for $k \in ( k_i, k_{i+1} ]$ as well.
If $ [k_i, k_{i+1}) \subset S_{bad}$ then from (\ref{phibd}) we have
\[
\| \phi ( j ) \| \leq ( 1 + \frac{\gamma_3}{\eps^{1/2}})  | d(j-1) | + | r(j-1) | , \;
j = k_i+1,..., k_{i+1} ,
\]
which means that (\ref{nice}) holds for $k \in ( k_i, k_{i+1} ]$.
On the other hand, if
$[k_i, k_{i+1}) \subset S_{good}$, then
$k_i -1 \in S_{bad}$; from (\ref{phibd}) we have that
\[
\| \phi ( k_i) \|  \leq ( 1 + \frac{\gamma_3}{\eps^{1/2}})  | d( k_i-1 ) | + | r (k_i -1 ) |  .
\]
Using (\ref{good}) to analyse the behaviour on $[k_i , k_{i+1}]$, we have
\begin{eqnarray*}
\| \phi (k ) \|  & \leq & 
 \gamma_4 \lambda^{k-k_i} \| \phi ( k_i  ) \| +  ( 2 \frac{ \| {\cal S} \|}{\delta} + 1 )
\gamma_4 \times \\
&&  \sum_{j=k_i}^{k-1}
\lambda^{k-1-j} ( | r(j) | + | d(j) | )  \\
&\leq &
\gamma_4 \lambda^{k-k_i } [ ( 1 + \frac{\gamma_3}{\eps^{1/2}})  | d( k_i-1 ) | + | r (k_i -1 ) | ] + \\
&& 
\gamma_4 ( 2 \frac{ \| {\cal S} \|}{\delta} + 1 ) \sum_{j=k_i}^{k-1}
\lambda^{k-1-j} ( | r(j) | + | d (j) | )   \\
& \leq & [ \gamma_4 ( 1 + \frac{\gamma_3}{\eps^{1/2}}) + \gamma_4 ( 2 \frac{ \| {\cal S} \|}{\delta} + 1 ) ] \times \\
&& 
\sum_{j=k_i-1}^{k-1}
\lambda^{k-1-j}  ( | r(j) |  + | d(j) | ) \\
& \leq & \gamma_5 
\sum_{j= t_0 } ^{k-1}  \lambda^{k-1-j} ( | d(j) | + | r(j) | )  , \\
& &
\;\;\; k = k_i+1 ,..., k_{i+1} ,
\end{eqnarray*}
as desired.
$\Box$

This completes the proof.

\noindent
$\Box$

\section{Tolerance to Time-Variations}

The linear-like bound proven in Theorem 1 can be leveraged to prove that
the same behaviour will result even in the presence of slow time-variations with
occasional jumps.
So suppose that the actual plant model is
\be
y(t+1) = \phi(t)^T \theta^* (t) +d(t) , \; \phi (t_0) = \phi_0 , 
\label{newplant}
\ee
with $\theta^* (t) \in {\cal S}$ for all $t \in \R$. We adopt a common
model of acceptable time-variations used in adaptive control:
with $c_0 \geq 0$ and $\eps >0$, we let
$s( {\cal S} , c_0,  \eps )$ denote the subset of
$\linf( \R^{2n})$ whose elements $\theta^*$ satisfy
$\theta^* (t) \in {\cal S}$ for every $t \in \Z$ as well as
\be
\sum_{t=t_1}^{t_2-1} \| \theta^* (t+1) - \theta^* (t) \| \leq c_0 + \eps ( t_2 - t_1 ) , \;
t_2 > t_1 
\label{tsbound}
\ee
for every $t_1 \in \Z$.
We will now show that, for every $c_0 \geq 0$, the approach tolerates time-varying parameters in
$s( {\cal S} , c_0,  \eps )$ if $\eps$ is small enough.

\noindent
\framebox[85mm][l]{
\parbox{82mm}{
\begin{theorem}
For every
$\delta \in (0, \infty ]$,
$\lambda_1 \in ( \underline{\lambda} ,1)$ and $c_0 \geq 0$, there exists a $c_1>0$ and
$\eps >0$ so that for
every $t_0 \in \Z$, $\theta_0 \in {\cal S}$,
$\theta^* \in s( {\cal S} , c_0, \eps )$, $\phi_0 \in \R^{2n}$,
and $y^*, d \in \ell_{\infty}$, 
when the adaptive controller 
(\ref{estimator1}), (\ref{estimator2})
and
(\ref{conguy}) is applied
to the time-varying plant (\ref{newplant}),
the
following holds:
\[
\| \phi (k) \| \leq c_1 \lambda_1^{k- t_0} \| \phi_0 \| +
\sum_{j=t_0}^{k-1}  c_1 \lambda_1^{k-1-j} {( | r(j) | + | d(j)| )} , \;\;
\]
\[
\;\;\;\;\;\;\;\;
k \geq t_0 .
\]
\end{theorem}
}}

\noindent
{\bf Proof:}

Fix $\delta \in ( 0, \infty ]$, $\lambda_1 \in (  \underline{\lambda} ,1)$,
$\lambda \in ( \underline{\lambda} , \lambda_1 )$ and $c_0 > 0$. 
Let $t_0 \in \Z$, $\theta_0 \in {\cal S}$, $\phi_0 \in \R^{2n}$,
and $y^*, d \in \ell_{\infty}$
be arbitrary. 
With $m \in \N$, we will consider $\phi (t)$ on intervals of the form
$[t_0 + im , t_0 + (i+1)m ]$; we will be analysing these intervals in
groups of $m $ (to be chosen shortly); we set $\eps = \frac{c_0}{m^2}$,
and let $\theta^* \in s( {\cal S} , c_0, \eps )$ be arbitrary.

First of all, for $i \in \Z^+$ we can 
rewrite the plant equation as
\begin{eqnarray}
y(t+1) &=& \phi (t)^T \theta^* (t_0 + im) + d(t) + \nonumber \\
&&
\underbrace{ \phi (t)^T [\theta^* ( t) - \theta^* (t_0 + im) ]}_{=: \tilde{n} (t)} , \;  
\nonumber \\
&& \;\;\;\; t \in [ t_0 + im , t_0+ (i+1)m ] .
\label{nplant}
\end{eqnarray}
Theorem 1 applied to (\ref{nplant}) says that
there exists a constant $c>0$ so that
\[
\| \phi (t) \|
\leq
c \lambda ^{ t - t_0 - im } \| \phi (t_0 + im ) \| +
\]
\[
\;\;\;
\sum_{j=t_0 + im }^{t-1} c  \lambda^{t-1-j} ( | r (j) | + | d(j) |  + | \tilde{n} (j) | ) , \;
\]
\[
\;\;\;\;\;\;\;\;\;\;
 t \in [ t_0 + im , t_0+ (i+1)m ] .
\]
The above is a difference inequality associated with
a first order system; using this observation together with the fact that
$c \geq 1$, we see that if we define
\[
\psi (t+1) = \lambda \psi (t) + | r(t) | + | d(t) | + | \tilde{n} (t) | , \;
\]
\[
\;\;\;\;\;
t \in [ t_0 + im , t_0+ (i+1)m-1 ] ,
\]
with $\psi (t_0 + im) = \| \phi (t_0 + im ) \|$, then
\[
\| \phi (t) \| \leq c \psi (t) , \;\; t \in [ t_0 + im , t_0+ (i+1)m ] .
\]
Now we analyse this equation for $i =0,1,..., m-1$.

\noindent
{\bf Case 1:} $|\tilde{n} (t) | \leq 
\frac{1}{2c} ( \lambda_1 - \lambda ) \| \phi ( t ) \| $ for 
all $t \in [ t_0 + im , t_0+ (i+1)m ]$.

In this case 
\begin{eqnarray*}
\psi (t+1) & \leq & \lambda \psi (t) + |r(t)| + | d(t) | + | \tilde{n} (t) | \\
%& \leq &  \lambda \psi (t) + |r(t)| + | d(t) |  + \frac{1}{2c} ( \lambda_1 - \lambda ) 
%\| \phi ( t  ) \| \\
& \leq &  \lambda \psi (t) +  |r(t)| + | d(t) |  + \frac{1}{2c} ( \lambda_1 - \lambda  )
c \psi (t) \\
& \leq & ( \frac {\lambda + \lambda_1 }{2} ) \psi (t) + | r(t) | +  | d(t) | , \\
&& \;\;\;
\; t \in [ t_0 + im , t_0+ (i+1)m ] ,
\end{eqnarray*}
which means that
\[
| \psi (t) | \leq (\frac{\lambda + \lambda_1 }{2} )^{t-t_0-im}
| \psi ( t_0 + im ) | +
\]
\[
\;\;
\sum_{j=t_0+im}^{t-1} (\frac{\lambda + \lambda_1 }{2} )^{t-1-j} ( | r (j) | + | d(j) | ) ,
\]
\[
\;\;\;\;\; t=t_0+im,..., t_0+ (i+1)m .
\]
This, in turn, implies that
\[
\| \phi ( t_0+ (i+1)m ) \|  \leq  c (\frac{\lambda + \lambda_1 }{2} )^m
\| \phi ( t_0 + im) \| + \nonumber 
\]
\be
  \sum_{j=t_0 + im}^{t_0 + (i+1)m-1}
c (\frac{\lambda + \lambda_1 }{2} )^{t_0 + (i+1)m-1-j} ( | r(j) | + | d(j) | ) .
\label{newgood}
\ee

\vspace{0.3cm}
\noindent
{\bf Case 2:} $|\tilde{n} (t) | >
\frac{1}{2c} ( \lambda_1 - \lambda ) \| \phi ( t ) \| $ for          
some
$t \in [ t_0 + im , t_0+ (i+1)m ]$.

Since $\theta^* (t) \in {\cal S}$ for $t \geq t_0$, we see
\[
|\tilde{n} (t) |  \leq 2 \| {\cal  S} \| \times \| \phi ( t ) \| , \;
t \in [ t_0 + im , t_0+ (i+1)m ] .
\]
This means that
\begin{eqnarray*}
\psi (t+1) & \leq & \lambda \psi (t) + | r(t)|+ | d(t) | + | \tilde{n} (t) | \\
%& \leq &  \lambda \psi (t) + | r(t)|+ | d(t) |  + 2 \| {\cal  S} \| \times \| \phi ( t -1 ) \| \\
& \leq &  \lambda \psi (t) + | r(t)|+ | d(t) |  + 2 \| {\cal  S} \| c
\psi (t) \\
& \leq & \underbrace{( 1+2c \| {\cal S} \|  )}_{=: \gamma_1}
\psi (t) + |r(t)| + | d(t) | , \\
&& \;\;
\; t \in [ t_0 + im , t_0+ (i+1)m ] ,
\end{eqnarray*}
which means that
\[
| \psi (t) | \leq \gamma_1^{t-t_0-im}
\| \psi ( t_0 + im ) \| +
\]
\[
\;\;
\sum_{j=t_0+im}^{t-1} \gamma_1^{t-j-1} (|r(j)| + | d (j) |) ,
\; t=t_0+im,..., t_0+ (i+1)m .
\]
This, in turn, implies that
\begin{eqnarray}
&& \| \phi ( t_0+ (i+1)m ) \| 
 \leq  c \gamma_1 ^m
\| \phi ( t_0 + im) \| + \nonumber \\
&& \;  c \sum_{j=t_0 + im}^{t_0 + (i+1)m-1}
( \gamma_1 )^{t_0 + (i+1)m-j-1} ( | r(j) | + | d(j) | ) \nonumber \\
& & \leq 
c \gamma_1^m \| \phi ( t_0 + im) \| +
 c ( \frac{ 2 \gamma_1}{ \lambda + \lambda_1 } )^m \times \nonumber 
\end{eqnarray}
\be
 \sum_{j=t_0 + im}^{t_0 + (i+1)m-1} ( \frac{ \lambda + \lambda_1}{2} )^
{t_0 + (i+1)m -j-1 } ( | r(j) | + | d(j) | ) .
\label{bad}
\ee

On the interval $[t_0 , t_0 + m^2]$ there are $m$ sub-intervals of length $m$;
furthermore, because of the choice of $\eps$ we have that
\[
\sum_{j=t_0}^{t_0 + m^2-1} \| \theta^* (j+1) - \theta^* (j) \|
\leq c_0 + \eps m^2  \leq 2 c_0 .
\]
A simple calculation reveals that there are at most
$N_1 := \frac{4 c_0 c}{ \lambda_1 - \lambda}$
sub-intervals
which fall into the categorory of Case 2, with the remaining number
falling into the category of Case 1.
Henceforth we assume that $m > N_1$.
If we use (\ref{newgood}) and (\ref{bad})
to analyse the behaviour of the closed-loop system on the interval
$[t_0 , t_0 + m^2]$, we end up with a crude bound of
\begin{eqnarray}
\| \phi (t_0 + m^2) \| & \leq &
c^{m} \gamma_1^{N_1 m}
( \frac{ \lambda_1 + \lambda}{2} )^{m (m- N_1)} \| \phi ( t_0 ) \| + \nonumber 
\end{eqnarray}
\[
  ( \frac{ 2 \gamma_1 }{ \lambda + \lambda_1 } )^{ m} 
( c \gamma_1^m)^m ( \frac{2}{ \lambda+ \lambda_1 } )^{m^2} \times
\]
\be
\sum_{j=t_0}^{t_0 +m^2-1} (\frac{ \lambda_1 + \lambda}{2} )^{t_0 +m^2-j-1}
( | r(j) | + | d(j) | ) .
\label{phibd2}
\ee
At this point we would like to choose $m$ so that
\[
c^{m} \gamma_1^{N_1 m} ( \frac{ \lambda_1 + \lambda}{2} )^{m^2 - m N_1} \leq \lambda_1^{m^2}  \; 
\]
\[
\Leftrightarrow \;
c^{m} \gamma_1^{N_1 m} ( \frac{2}{ \lambda+ \lambda_1 })^{m N_1} \leq
( \frac{2 \lambda_1}{\lambda_1 + \lambda } )^{m^2};
\]
notice that $\frac{2 \lambda_1}{\lambda_1 + \lambda } > 1$, so if
we take the log of both sides, we see that we need
\begin{eqnarray*}
&& m \ln (c) + N_1 m \ln ( \gamma_1 ) + N_1 m \ln ( \frac{2}{ \lambda + \lambda_1 } )  \\
& \leq &
m^2 \ln ( \frac{2 \lambda_1}{\lambda_1 + \lambda } ) ,
\end{eqnarray*}
which will clearly be the case for large enough $m$, so at this point
we choose such an $m$.
It follows from (\ref{phibd2}) that there exists a constant
$\gamma_2$ so that
\[
\| \phi (t_0 + m^2) \|  \leq \lambda_1^{m^2} \| \phi (t_0) \| +
\]
\[
\;\;\;
\gamma_2 \sum_{j=t_0}^{t_0 +m^2-1} \lambda_1^{t_0 +m^2-j-1 } ( | r(j) | +  | d(j) | ) .
\]
Indeed, by time-invariance of the closed-loop system we see that
\[
\| \phi (\bar{t} + m^2) \|  \leq \lambda_1^{m^2} \| \phi ( \bar{t}) \| +
\]
\[
\;\;
\gamma_2 \sum_{j=\bar{t}}^{\bar{t} +m^2-1} \lambda_1^{\bar{t} +m^2-j-1 } ( | r(j) | +  | d(j) | ) , \; \bar{t} \geq t_0 .
\]
Solving iteratively yields
\[
\| \phi (t_0 + i m^2) \|  \leq \lambda_1^{im^2} \| \phi (t_0) \| +
\]
\be
\gamma_2 \sum_{j=t_0}^{t_0 +im^2-1} \lambda_1^{t_0 +im^2-j-1 } ( | r(j)| + | d(j) | ) , \; i \geq 0 .
\label{sampled}
\ee
We now combine this bound with the bounds which hold on the good intervals
(\ref{newgood}) and the bad intervals (\ref{bad}), and conclude that there
exists a constant $\gamma_3$ so that
\[
\| \phi (t ) \| \leq \gamma_3 \lambda_1^{t- t_0} \| \phi (t_0 ) \| +
\]
\[
\;\;
 \gamma_3 \sum_{j=t_0}^{t-1} \lambda_1^{t-j-1} (| r(j) | + | d(j) | ) , \;
t \geq t_0 ,
\]
as desired.
$\Box$

\section{Tolerance to Unmodelled Dynamics}

Due to the linear-like bounds
proven in Theorems 1 and 2, we can use the Small Gain Theorem
to good effect to prove the tolerance of the closed-loop system to unmodelled dynamics.
However, since the controller, and therefore the closed-loop system, is nonlinear, handling 
initial conditions is more
subtle: in the linear-time invariant case we can separate out the effect of
initial conditions from that of 
the forcing functions ($r$ and $d$), but in our situation they are inter-twined.
We proceed by looking at two cases - with and without 
initial conditions.
In all of the cases we consider the time-varying plant (\ref{newplant}) 
with $d_{\Delta} (t)$ added to represent the effect of unmodelled dynamics:
\be
y(t+1) = \phi(t)^T \theta^* (t) +d(t) + d_{\Delta} (t) , \; \phi (t_0) = \phi_0 .
\label{mod_dude}
\ee
To proceed,
fix $\delta \in (0, \infty ]$, $\lambda_1 \in ( \underline{\lambda} ,1)$ and $c_0 \geq 0$; 
from Theorem 2
there exists a 
$c_1>0$ and
$\eps >0$ so that for
every $t_0 \in \Z$, $\phi_0 \in \R^{2n}$, $\theta_0 \in {\cal S}$,
$y^*, d \in \ell_{\infty}$,  and $\theta^* \in s( {\cal S} , c_0, \eps )$,
when the adaptive controller
(\ref{estimator1}), (\ref{estimator2})
and
(\ref{conguy}) is applied
to the time-varying plant (\ref{mod_dude}),
the
following bound holds:
\begin{eqnarray}
\| \phi (k) \| & \leq & c_1 \lambda_1^{k- t_0} \| \phi_0 \| + \nonumber \\
&& \sum_{j=t_0}^{k-1}  c_1 \lambda_1^{k-1-j} {( | r(j) | + | d(j)| + | d_{\Delta} (j) | )} , \;\; \nonumber \\
&& \;\;\;\;
k \geq t_0 .
\label{goodbound}
\end{eqnarray}

\subsection{Zero Initial Conditions}

In this case we assume that $\phi (t) = 0$ for $t \leq t_0$; we derive a
bound on the closed-loop system behavour in the presence of unmodelled dynamics.
Suppose that the unmodelled dynamics is
of the form
$d_{\Delta} (t) = ( \Delta \phi )(t)$ with
$\Delta : \linf ( \R^{2n})  \rightarrow \linf (\R^{2n})$ a (possibly
nonlinear time-varying) causal map 
with a finite gain of $\| \Delta \|$.
It is easy to prove that if
%\[
%\frac{c_1}{1- \lambda_1} \| \Delta \|  < 1 \; \Leftrightarrow \;
$\| \Delta \| < \frac{1 - \lambda_1}{c_1  } $,
then
\[
\| \phi (k)  \| \leq \frac{c_1}{1 - \lambda_1 - c_1 \| \Delta \|  }
( \sup_{t \geq t_0} \| r (t)  \|+
\sup_{t \geq t_0} \| d(t) \|  ) , 
\]
\[
\;\;\;\;\; k \geq t_0 ,
\]
i.e., a form of closed-loop stability is attained.
Following the approach of Remark 5, we could also analyse the closed-loop
system using $l_p$-norms with $1 \leq p < \infty$.

\subsection{Non-Zero Initial Conditions}

Now we allow unmodelled LTI dynamics with non-zero initial conditions, and we
develop convolution-like bounds on the closed-loop system.
To this end suppose that
the unmodelled dynamics are of the form
\be
d_{\Delta} (t) :=
\sum_{j=0}^{\infty} \Delta_j  \phi(t-j) ,
\label{mod1}
\ee
with $\Delta_j  \in \R^{1 \times 2n }$;
the corresponding transfer function is $\Delta (z^{-1}) := \sum_{j=0}^{\infty}
\Delta_j   z^{-j} $.
It is easy to see that this model subsumes the classical additive uncertainty,
multiplicative uncertainty, and uncertainty in a coprime
factorization, which is common in the robust control literature, e.g. see
\cite{zhou_book}, with the only
constraint being that the perturbations correspond to strictly causal terms.
In order to obtain linear-like bounds on the closed-loop behaviour,
we need to impose more constraints on $\Delta (z)$ than in the previous sub-section:
after all,
if $\Delta (z^{-1}) = \Delta_p  z^{-p}$, it is clear that $\| \Delta \| = 
\| \Delta_p  \|$
for all $p$, but the effect on the closed-loop system varies greatly - a large value of
$p$ allows the behaviour in the far past to affect the present.
To this end, with $\mu >0$ and $\beta \in (0,1)$, 
we shall restrict $\Delta (z^{-1})$ to a set of the form
\[
{\cal B} ( \mu , \beta ) :=
\{  \sum_{j=0}^{\infty}
\Delta_j   z^{-j} : \Delta_j \in \R^{1 \times 2n} \mbox{ and }
\]
\[
 \;\;\;\;\;\;\;\; \| \Delta_j  \| \leq \mu \beta^j , \; j \geq 0 \} .
\]
It is easy to see that every transfer function in $ {\cal B} ( \mu , \beta ) $ is analytic
in $\{ z \in \C : | z |  > \beta \}$, so it has no poles
in that region.

Now we fix $\mu >0$ and $\beta \in (0,1)$ and let
$\Delta ( z^{-1}) $ belong to ${\cal B} ( \mu , \beta )$; the goal
is to analyse the closed-loop behaviour of (\ref{mod_dude}) for $t \geq t_0$
when $d_{\Delta}$ is given by (\ref{mod1}).
We first partition $\dd (t) $ into two parts - that which depends on
$\phi (t)$ for $t \geq t_0$ and that which depends on
$\phi (t)$ for $t < t_0$:
\begin{eqnarray*}
\dd (t) &=& 
\sum_{j=0}^{\infty} \Delta_j \phi(t-j) 
= \sum_{j= - \infty}^{t} \Delta_{ t-j} \phi (j)  \\
&=& \underbrace{\sum_{j=- \infty}^{t_0-1} \Delta_{ t-j} \phi(j)}_{=: \dd^- (t)} + 
\underbrace{\sum_{j=t_0}^{t} \Delta_{ t-j} \phi(j)}_{=: \dd^+ (t) } .
\end{eqnarray*}
It is clear that
\begin{eqnarray*}
\| \dd^+ (t) \| & \leq & \sum_{j=t_0}^{t} \mu \beta^{t-j} \| \phi(j) \| , \\
\| \dd^- (t) \| & \leq & \sum_{j=- \infty}^{t_0-1} \mu \beta^{t-j} \| \phi(j) \|  
\\
&=& \mu \beta^{t-t_0} \sum_{j=1}^{ \infty} \beta^j \| \phi ( t_0 - j ) \| , \;
t \geq t_0 .
\end{eqnarray*}
If $\phi (t)$ is bounded on $\{ t \in \Z: t < t_0 \}$ then
$\sum_{j=1}^{\infty}  \beta^{j} \| \phi( t_0-j) \|$ is finite, in which case
we see that $\dd^- (t)$ goes to zero exponentially fast; henceforth, we make the
reasonable assumption that this is the case.
It turns out that we can easily bound $\dd (t)$ with a difference equation. 
To this end, consider
\be
{m} (t+1) = \beta   {m} (t) + \beta \| \phi (t) \| , \; t \geq t_0 ,
\label{meq}
\ee
with ${m} (t_0) = m_0 := \sum_{j=1}^{\infty}  \beta^{j} \| \phi( t_0-j) \|$; it is straight-forward
to prove that
\be
| \dd (t) | \leq | \dd^+ (t) | + | \dd^- (t)| \leq  \mu  {m} (t) + \mu \| \phi (t) \|  , \; t \geq t_0 .
\label{meq2}
\ee
This model of unmodelled dynamics is similar to that
used in the adaptive control literature, e.g.
see \cite{kreiss}.

\noindent
\framebox[85mm][l]{
\parbox{82mm}{
\begin{theorem}
For every $\beta \in (0,1)$ and
$\lambda_2 \in ( \max\{\lambda_1 , \beta \} , 1 )$, there 
exist $\bar{\mu} > 0$ and $c_2 > 0$ so that 
for every $t_0 \in \Z$, $\phi_0 \in \R^{2n}$, $m_0 \in \R$, $\theta_0 \in {\cal S}$,
$y^* , d \in \linf$, $\theta^* \in s ( {\cal S}, c_0 , \eps )$ and
$\mu \in (0, \bar{\mu} )$,
when the adaptive controller
(\ref{estimator1}), (\ref{estimator2})
and
(\ref{conguy}) is applied
to the time-varying plant (\ref{mod_dude}) with $d_{\Delta}$ satisfying
(\ref{meq}) and (\ref{meq2}),
%(\ref{mod1}), $\Delta (z^{-1}) = \sum_{j=0}^{\infty} \Delta_j z^{-j}$
%lying in ${\cal B} ( \mu , \beta )$,
the following bound holds:
\[
\| \phi (k) \| \leq c_2 \lambda_2^{k - t_0} ( \| \phi_0 \| + | m_0 | ) +
\]
\[
\;\;\;
\sum_{j = t_0}^{k-1} c_2 \lambda_2^{k - 1-  j}  ( | d(j) | + | r(j) | ) ,
 \;\; k \geq t_0 .
\]
\end{theorem}
}}

\noindent
{\bf Proof:}

Fix $\beta \in (0,1)$ and $\lambda_2 \in (  \max\{\lambda_1 , \beta \} , 1 )$.
The first step is to convert difference inequalities to difference
equations. To this end, consider the 
difference equation
\[
\tphi (t+1) = \lambda_1 \tphi (t) + c_1 | r(t) | + c_1 | d(t) | + c_1
 \mu \tilde{m} (t) +
\]
\be
\;\;\;  c_1 \mu \tilde{\phi }(t) , \;
\tphi (t_0) = c_1 \| \phi (t_0) \| ,
\label{eqn1}
\ee
together with the difference equation based on (\ref{meq}):
\be
\tilde{m} (t+1) = \beta   \tilde{m} (t) + \beta  \tphi (t)  , \; 
\tilde{m} (t_0) = | {m}_0 | .
\label{eqn2}
\ee
It is easy to use induction together with (\ref{goodbound}), (\ref{meq}), and (\ref{meq2}) 
to prove that
\be
 \| \phi (t) \| \leq  \tphi (t) , 
\;\; | m(t) | \leq \tilde{m} (t) , \; t \geq t_0 . 
\label{phimbd}
\ee
If we combine the difference equations (\ref{eqn1}) with 
(\ref{eqn2}), we end up with 
\[
\tbo{ \tphi (t+1)}{\tilde{m}(t+1)} =
\underbrace{\tbt{ \lambda_1 +  c_1\mu}{ c_1 \mu}{ \beta } {\beta}}_{A_{cl} ( \mu )}
\tbo{ \tphi (t)}{\tilde{m}(t)} +
\]
\be
\;\;\;
 \tbo{c_1 }{0} ( | d(t) |  + | r(t) | ), \;
t \geq t_0 .
\label{phit}
\ee
Now we see that $A_{cl} ( \mu ) \rightarrow \tbt{ \lambda_1}{ 0}{\beta }{\beta }$
as $\mu \rightarrow 0$, and this matrix has eigenvalues of
$ \{ \lambda_1 , \beta  \}$.
Now choose $\bar{\mu} > 0$ so that all eigenvalues are less than
$( \frac{\lambda_2}{2} + \frac{1}{2}\max \{ \lambda_1 , \beta \} )   $ in magnitude for 
$\mu \in (0 , \bar{\mu} ]$, and define $\eps := 
\frac{\lambda_2}{2} - \frac{1}{2}\max \{ \lambda_1 , \beta \}$. 
Using the proof technique of Desoer in \cite{desoer}, we can conclude that
for $\mu \in  (0 , \bar{\mu} ]$, we have
\[
\| A_{cl} ( \mu ) ^k \| \leq \underbrace{\left( \frac{3 + 2\beta + 2 c_1 \bar{\mu}}{ \eps^2} \right)}_{=: \gamma_1} 
 \lambda_2^k , \;\; k \geq 0 ;
\]
if we use this in (\ref{phit}) 
and then apply the bounds in (\ref{phimbd}), it follows that
%and then 
%which means that
%\[
%\left  \| \tbo{  \tphi (t)}{\tilde{m}(t)} \right \|
%\leq \gamma_1 \lambda_2^{t- t_0} \left \| \tbo{  \tphi (t_0)}{\tilde{m}(t_0)} \right \| +
%
%\sum_{j=t_0}^{t-1} \gamma \lambda_2^{t-j-1}(  | d (j) | + | r(j)| )  , \; t \geq t_0 ;
%\]
%using the bounds given above, it follows that
\[
\| \phi (k) \| \leq c_1 \gamma_1 \lambda_2 ^{ k - t_0} ( \| \phi_0 \| + | m_0 | ) +
\]
\[
\;\;\;
\sum_{j=t_0}^{k-1} c_1 \gamma_1 \lambda_2^{k-1-j} ( | d (j) | + | r(j) | )  , \;\; k \geq t_0,
\]
as desired.
$\Box$

\section{Step Tracking}

If the plant is non-minimum phase, it is not possible track an arbitrary
bounded reference
signal using a bounded control signal.
However, as long as the plant does not have a zero at $z=1$, 
it is possible to modify the controller design procedure to achieve asymptotic step tracking
if there is no noise/disturbance. 
So at this point assume that the corresponding plant polynomial $B( z^{-1} )$
has no zero at $z=1$ for any plant model $\theta^* \in {\cal S}$.
To proceed, we use
the standard trick from the literature, e.g.
see \cite{goodwinsin}: we still estimate $A(z^{-1} )$ and $B( z^{-1}) $
as before, but we now design the control law slightly differently. To this end, we first
define 
\[
\tilde{A} (t, z^{-1} ) := (1- z^{-1} ) \hat{A} ( t,  z^{-1} ) ,
\]
and then
let $A^* ( z^{-1} )$ be a $2 (n+1)^{th}$ monic polynomial (rather than a $2n^{th}$ one)
of the form
\[
A^* (z^{-1} ) = 1 + a_1^* z^{-1} + \cdots + a_{2n+2}^* z^{-2n-2}  
\]
so that $z^{2(n+1)} A^* (z^{-1} )$ has all of its zeros in $\D^o$.
Next, we choose two polynomial
\[
{\tilde{L}} (t, z^{-1} ) = 1 + \tilde{l}_1 (t) z^{-1} + \cdots + \tilde{l}_{n+1} (t)
z^{-n-1}  
\]
and
\[
\hat{P} (t,z^{-1} ) = \hat{p}_1 (t) z^{-1} + \cdots + \hat{p}_{n+1} (t) z^{-n-1}
\]
which satisfy the equation
\be
\tilde{A} ( t , z^{-1} ) \tilde{L} (t, z^{-1} ) +
\hat{B} ( t , z^{-1} ) \hat{P} (t,z^{-1} ) = A^* ( z^{-1} ) ;
%\label{poleplace}
\ee
since
$\tilde{A} ( t , z^{-1} ) $ and
$\hat{B} ( t , z^{-1} )$ are coprime, there 
exist {\bf unique}
$\tilde{L} (t, z^{-1} )$ and $\hat{P} (t,z^{-1} )$ which satisfy this
equation.
We now define 
\[
\hat{L} (t, z^{-1} ) =  
1 + \hat{l}_1 (t) z^{-1} + \cdots + \hat{l}_{n+2} (t) z^{-n-2} 
\]
\[
\;\; := 
(1-z^{-1} ) \tilde{L} (t, z^{-1} ) ;
\]
at time $t$ we choose $u(t)$ so that
\[
  u(t) =  - \hat{l}_1  (t-1) u(t-1) - \cdots - 
\]
\[
 \hat{l}_{n+2} (t-1) u(t-n-2)
- \hat{p}_1 (t-1) [ y(t-1)- y^* (t-1) ] - \cdots  
\]
\[
  - \hat{p}_{n+1} (t-1) [ y(t-n-1)- y^* (t-n-1) ]  .
\]
We can use a modified version
of the argument used in the proof of Theorem 1 to conclude that a similar type of
result holds here; 
we can also prove that
asymptotic step tracking will be attained if the noise is zero
and the reference signal $y^*$ is constant.
The details are omitted due to space considerations.

\section{A Simulation Example}

Here we provide an example to illustrate the benefit of the proposed adaptive controller.
To this end, consider the second order plant
\[
y(t+1) = - a_1(t) y(t) - a_2 (t) y(t-1)+ 
\]
\[
\;\;\;\;b_1 (t) u(t) + b_2 (t) u (t-1) + d(t)
\]
with $a_1  (t)\in [0 , 2]$, $a_2 (t) \in [1,  3 ]$, $b_1 (t) \in [0,1]$, and $b_2 (t) \in [-5, -2]$.
So every admissible model is unstable and non-minimum phase, which makes this a 
challenging plant to control. We set $\delta = \infty$.

\subsection{Stability}

In this sub-section we consider the problem of stability only - we set $y^* = 0$.
First we compare the ideal algorithm (\ref{orig2})-(\ref{est2})
(with projection onto ${\cal S}$)
with the classical one (\ref{new})
(suitably modified to have projection onto ${\cal S}$); in both cases we couple
the estimator with the adaptive pole placement controller (\ref{conguy}) where we place all
closed-loop poles at zero.
In the case of the classical estimator (\ref{new}) we arbitrarily set
$\alpha = \beta =1$.
Suppose that
the actual value of $(a_1,a_2,b_1,b_2) $ is $(2,3,1,-2)$ and the initial estimate is set
to the midpoint of the interval.
In the first simulation we
set $y(0) = y(-1) =0.01$ and $u(-1)=0$ and set the noise $d(t)$ to zero -
see the top plot of Figure 1.
In the second simulation we set
$y(0)= y(-1)=u(-1)=0$ and the noise to
$d(t) = 0.01*\sin (5 t)$ - see the bottom
plot of Figure 1.
In both cases the controller based on the
ideal algorithm (\ref{orig2})-(\ref{est2}) is clearly superior to the one based on the
revised classical
algorithm (\ref{new}).

\begin{figure}[htbp]
\begin{center}
\includegraphics[height=2.5in,width=3.5in]{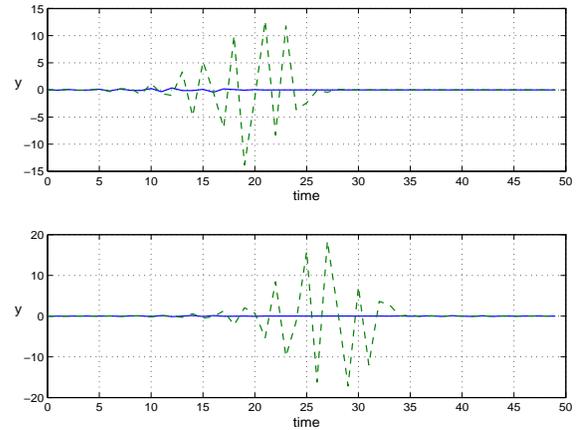}
        \caption{
A comparison of the ideal algorithm (solid) and the classical algorithm (dashed) with a non-zero initial condition and no noise (top plot) and a zero initial condition and noise (bottom plot).}
\end{center}
\end{figure}

Now we further examine the case of the proposed controller 
when it is applied to the time-varying plant with unmodelled
dynamics, a zero initial condition,
and a non-zero noise. 
More specifically, we
set
\[
a_1 (t) = 1 + \sin (.001t), \; a_2 (t) = 2 +  \cos(.001t) , \;
\]
\[
b_1 (t) = 0.5 + 0.5 \sin( .005t) , 
\;\;  b_2 (t) = -3.5 -1.5 \sin (.005t) , \;\;
\]
\[
d(t) = 0.01 \sin (5t) .
\]
For the unmodelled part of the plant we use a term of the form
discussed in Section VI.B:
\[
m(t+1) = 0.75 m(t) + 0.75 \| \phi (t) \| , \;\; m(0) = 0 ,
\]
\[
d_{
\Delta} (t) = \left \{ \begin{array}{ll}
0 & \;\; t=0,1,..., 4999 \\
0.025 m(t) + 0.025 \| \phi (t) \| & \;\; t \geq 5000.
\end{array} \right.
\]
We plot the result in Figure 2; we see that the parameter estimator approximately follows the
system parameters, and the effect of the noise is small on average,
even in the presence of unmodelled dynamics.

\begin{figure}[htbp]
\begin{center}
\includegraphics[height=3.0in,width=3.5in]{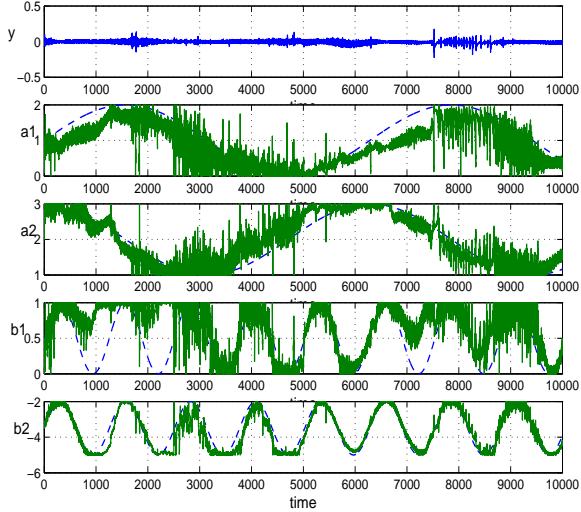}
        \caption{The system behaviour with time-varying parameters and unmodelled dynamics; the parameters are dashed and the estimates are solid.}
\end{center}
\end{figure}

\subsection{Step Tracking}

The plant in the previous sub-section 
has a large amount of uncertainty, as well as a wide range of unstable
poles and non-minimum phase zeros, which means that there are limits on
the quality of
the transient behaviour even if the parameters were fixed and known.
Hence, to illustrate the tracking ability we look at a sub-class of 
systems: one with $a_1$ and $b_1$ as before, namely
$a_1  (t)\in [0 , 2]$ and $b_1 (t) \in [0,1]$, but now with
$a_2 = 1$ and $b_2 = -3.5$.
With fixed parameters the corresponding system is still unstable and
non-minimum phase.

We simulate the closed-loop pole placement step tracking controller 
of Section VII with
a zero initial condition,
initial parameter estimates at the midpoints of the admissible intervals, and
with time-varying parameters:
\[
a_1 (t) = 1 + \sin (.002t), 
b_1 (t) = 0.5 + 0.5 \cos( .005t) , 
\]
with a non-zero disturbance:
\[
d(t) = \left \{ \begin{array}{ll}
 0.01 \sin (5t) & \;\; t =0,1,..., 2499 \\
0.05 \sin (5t) & \;\; t = 2500,..., 4999,
\end{array} \right.
\]
and a square wave reference signal of
$y^*(t) = \mbox{sgn} [ sin (0.01t)]$.
We plot the result in Figure 3; we see that 
the parameter estimates crudely
follows the
system parameters, with less accuracy than in the previous sub-section,
partly due to the fact that the constant setpoint dominates the estimation
process and leads to higher inaccuracy. As a result,
$y(t)$ does a good job of following
$y^*$ on average, but with the occasional flurry of activity
when the parameter estimates are
highly inaccurate.
When the noise is increased five-fold at $k=2500$, the behaviour degrades
only slightly.

\begin{figure}[htbp]
\begin{center}
\includegraphics[height=2.5in,width=3.5in]{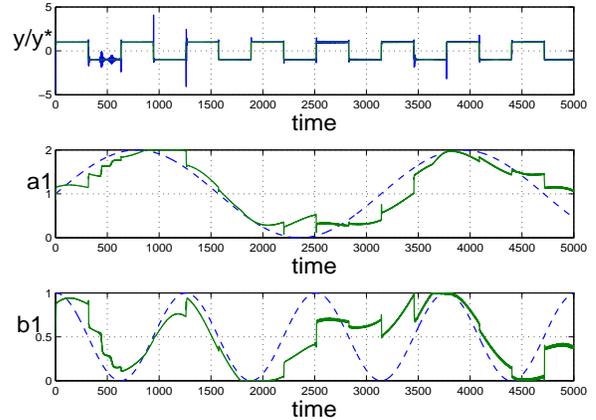}
        \caption{The pole placement tracking controller with time-varying parameters and small noise; the parameters are dashed and the estimates are solid.}
\end{center}
\end{figure}

\section{Summary and Conclusions}

Here we show that
if the original, ideal, projection algorithm is used in the
estimation process (subject to the assumption that the plant parameters lie in a convex,
compact set), then  the corresponding
pole placement adaptive controller
guarantees linear-like convolution bounds on the closed
loop behaviour, which confers exponential stability and a bounded
noise gain (in every $p$-norm with $1 \leq p \leq \infty$), unlike almost all other parameter adaptive controllers.
This can be leveraged to prove tolerance to unmodelled dynamics and
plant parameter variation.
We emphasize that there is no persistent excitation requirement of any sort;
the improved performance arises from the vigilant nature of the
the ideal parameter estimation algorithm.

As far as the author is aware, the linear-like convolution bound proven here
is a first in parameter adaptive control. It allows a modular approach
to be used in analysing time-varying parameters and
unmodelled dynamics. This approach avoids all of the fixes invented in the
1980s, such as signal normalization and
deadzones, used to deal with the lack of robustness to unmodelled
dynamics and time-varying parameters. 

We are presently working on extending the approach to
the model reference adaptive control setup. It will be interesting
to see if 
the convexity assumption can be removed by using
multi-estimators, i.e.
cover the the set
of admissible parameters by a finite number of convex sets, and then
use an estimator for each such set.
Extending the approach to the continuous-time setting may prove challenging,
since a direct application would yield a non-Lipschitz continuous estimator,
which brings with it mathematical solveability issues.

\section{Appendix}

\vspace{0.3cm}
\noindent
{\bf Proof of Proposition 1:}

Since projection does not make the parameter estimate worse,
it follows
from (\ref{estimator1}) that
\begin{eqnarray*}
&&
\| \hat{\theta} (t+1) - \hat{\theta} (t) \| 
 \leq  \| \check{\theta} (t+1) - \hat{\theta} (t) \|  \\
 & \leq & 
\|  \rhot 
\frac{ \phi (t)}{ \phi (t)^T  \phi (t)} e(t+1) \| \\
& \leq & 
 \rhot \frac{ |e(t+1)| }{ \|  \phi (t) \| } ,
\; t \geq t_0 .
\end{eqnarray*}
so the first inequality holds.
%which means that
%\be
%\sum_{j=t_0}^{t-1} \| \hat{\theta} (j+1) - \hat{\theta}  (j) \|
%\leq \sum_{j=t_0}^{t-1}  \rhoj \frac{ |e(j+1)| }{ \|  \phi (j) \| } .
%\label{errorsum}
%\ee

We now turn to energy analysis. 
We first define
$\tilde{\check{\theta}} (t):= \check{\theta} (t) - \theta^*$ and
%$V(t) :=  \tilde{\theta}(t)^T \tilde{\theta}(t) $ and
$\check{V} (t) := \tct (T)^T \tct (t) $.
Next, we 
subtract $\theta^*$ from each side of (\ref{estimator1}), yielding
\begin{eqnarray*}
\tct  (t+1) &=& \tilde{\theta} (t) 
+ \rhot 
\frac{ \phi (t)}{ \phi (t)^T  \phi (t)} \times \\
&&
 [ 
-\phi(t)^T  \tilde{\theta} (t)  + d(t) ] \\
&=& [ I - 
\underbrace{\rhot 
\frac{ \phi (t) \phi (t)^T}{ \phi (t)^T  \phi (t)}}_{=: W_1 (t)}  ]
\tilde{\theta} (t) + \\
&& 
\underbrace{ \rhot  \frac{ \phi (t)}{ \phi (t)^T  \phi (t)}}_{=: W_2 (t)} d(t) .
\end{eqnarray*}
Then
\begin{eqnarray*}
\check{V} (t+1) &=&
[ (I-W_1(t)) \ttt (t) + W_2 (t) d(t) ]^T \times \\
&&
[ (I-W_1(t)) \ttt (t) + W_2 (t) d(t) ] \nonumber \\
&=& \ttt (t)^T [ I - W_1(t)] [ I - W_1(t)] \ttt (t) + \\
&&   2 \ttt (t)^T 
[I - W_1(t)] W_2(t) d(t) + \\
&&  W_2(t)^T W_2(t) d(t)^2 .
\label{vcupdate}
\end{eqnarray*}
Now let us analyse the three terms on the RHS:
the fact that
$W_1(t)^2= W_1(t)$ allows us to simplify the first term;
the fact that $W_1 (t) W_2 (t) = W_2 (t)$ 
means that the second term is zero;
$W_2 (t)^T W_2 (t) = \rhot  \frac{1}{ \phi (t)^T  \phi (t)}$, which
simplifies the third term. We end up with
\[
\check{V} (t+1) 
\]
\begin{eqnarray*}
&=& 
\ttt (t)^T  [ I - W_1(t)]  \ttt (t)  + \rhot \frac{d(t)^2}{ \phi (t)^T  \phi (t)  } \nonumber \\
&=& V(t) - \rhot 
\frac{ [ \ttt (t)^T \phi (t)]^2}{ \phi (t)^T \phi (t)} +  \\
&& \;\;
 \rhot   \frac{d(t)^2}{ \phi (t)^T  \phi (t)} \nonumber \\
&=& V(t) + 
 \rhot   \frac{ d(t)^2 - [d (t)-e(t+1)]^2}
{ \phi (t)^T  \phi (t)} \label{vguy}  \\
& \leq & 
V(t) + \rhot   \frac{ - \frac{1}{2} e(t+1)^2 + 2 d(t)^2}
{\phi (t)^T  \phi (t)} .
\end{eqnarray*}
Since projection never makes the estimate worse, it follow that
\[
V(t+1) \leq V(t) + \rhot   \frac{ - \frac{1}{2} e(t+1)^2 + 2 d(t)^2}
{\phi (t)^T  \phi (t)} .
\]

\noindent
$\Box$

\noindent
{\bf Proof of Lemma 1:}
Fix $\delta \in ( 0 , \infty ]$ and $\sigma \in ( \underline{\lambda} , 1) $.
First of all, it is well known that the characteristic polynomial of
$\bar{A} (t)$ is exactly $z^{2n} A^* (z^{-1})$ 
for every $t \in \Z$.
Furthermore, it is well known that
the coefficients of
$\hat{L} (t, z^{-1} )$ and $\hat{P} (t,z^{-1} )$ are the solution of a linear
equation, and are analytic functions of
$\hat{\theta} (t) \in {\cal S}$.
Hence, there exists a constant $\gamma_1$ so that, for every set of initial
conditions, $y^* \in \linf$ and $d \in \linf$, we have
$\sup_{t \geq t_0} \| \bar{A} (t) \| \leq \gamma_1$.

To prove the first bound we now invoke the argument used in \cite{desoer},
who considered a more general time-varying situation but with more restrictions
on $\sigma$. By making a slight adjustment to the first part of the proof given there,
we can prove that with
$\gamma_2 := \sigma \frac{ (\sigma + \gamma_1 )^{2n-1}}{ ( \sigma - \underline{\lambda} )^{2n}}$,
then for every $t \geq t_0$ we have
$\| \bar{A}  (t) ^k \| \leq \gamma_2 \sigma^{k} , \;\; k \geq 0 $,
as desired.

Now we turn to the second bound. From Proposition 1 
and the
Cauchy-Schwarz inequality we obtain
\begin{eqnarray*}
&& \sum_{j=k}^{t-1} \| \hat{\theta}  (j+1) - \hat{\theta} (j) \|  \\
& \leq & 
\sum_{j=k}^{t-1} \rhoj \frac{| e(j+1) |}{\| \phi (j)\|}  \\
& \leq & [ \sum_{j=k}^{t-1} \rhoj \frac{e(j+1)^2}{\| \phi (j)\|^2} ]^{1/2}
( t-k)^{1/2} .
\end{eqnarray*}
Now notice that
\[
\| \bar{A} (t+1) - \bar{A} (t) \| \leq \| \hat{\theta} (t+1) - \hat{\theta}
(t) \| + 
\]
\[
\;\;\;\;\sum_{i=1}^n ( | \hat{l}_i (t+1) - \hat{l}_i (t) | + |
\hat{p}_i (t+1) - \hat{p}_i (t) | ) .
\]
The 
fact that 
the coefficients of 
$\hat{L} (t, z^{-1} )$ and $\hat{P} (t,z^{-1} )$ are analytic functions of
$\hat{\theta} (t) \in {\cal S}$
means that there exists a constant $\gamma_3 \geq 1$ so that
\[
\sum_{j=k}^{t-1} \| \bar{A} (j+1) - \bar{A} (j) \| \leq 
\gamma_3 \sum_{j=k}^{t-1} \| \hat{\theta}  (j+1) - \hat{\theta} (j) \|,
\]
so we conclude that the second bound holds as well.
$\Box$.


\begin{thebibliography}{999}

\bibitem{desoer}
C.A. Desoer,
``Slowly Varying Discrete Time System $x_{t+1} = A_t x_t$",
{\em Electronic Letters},
vol. 6, no. 11, pp. 339 - 340, May 1970.


\bibitem{morse1978}
A. Feuer and A.S. Morse,
``Adaptive Control of Single-input, Single-output Linear Systems",
{\em IEEE Transactions on Automatic Control}, vol. 23, No. 4, pp. 557-569, 1978.


\bibitem{barmish}
M. Fu and B.R. Barmish,
``Adaptive Stabilization of Linear Systems Via
Switching Control'',
{\em IEEE Transactions on Automatic Control},
vol. AC-31, pp. 1097-1103, Dec. 1986.

\bibitem{goodwin1980}
G.C. Goodwin,  P.J. Ramadge, and P.E.  Caines,
  ``Discrete Time Multivariable Control'',
  {\it IEEE Transactions on Automatic Control},
vol. AC--25, pp. 449--456, 1980.

\bibitem{goodwinsin}
G.C. Goodwin and K.S. Sin,
``Adaptive Filtering Prediction and Control, Prentice Hall,
Englewood Cliffs, New Jersey, USA, 1984.

\bibitem{liberzon1}
J. P. Hespanha, D. Liberzon, and A. S. Morse,
``Hysteresis-based switching algorithms for supervisory control of
uncertain systems",  {\em Automatica}, vol. 39, pp. 263-272, Feb 2003.


\bibitem{liberzon2}
J. P. Hespanha,
D. Liberzon, and A. S. Morse,
``Overcoming the limitations of adaptive control by means of
logic-based switching",
{\em Systems and Control Letters}, vol. 49, no. 1, pp. 49-65, May 2003.

\bibitem{Ioa86}
P.A. Ioannou and K.S. Tsakalis,
``A Robust Direct Adaptive Controller",
{\em IEEE Transactions on Automatic Control},
{\bf vol. AC-31}, no. 11, pp. 1033 -- 1043, 1986.


\bibitem{hanfu}
Y. Li and H-F Chen,
``Robust Adaptive Pole Placement for Linear Time-varying Systems",
{\em IEEE Transactions on Automatic Control},
vol. AC-41, pp. 714 -- 719, 1996.


\bibitem{kreiss2}
G. Kreisselmeier,
``Adaptive Control of a  Class of Slowly Time-varying Plants",
{\em Systems and Control Letters}, vol. 8, pp. 97 -- 103, 1986.

\bibitem{kreiss}
G. Kreisselmeier and B.D.O. Anderson,
``Robust Model Reference Adaptive Control",
{\em IEEE Transactions on Automatic Control}, {AC-31},
pp. 127 -- 133, 1986.

%\bibitem{lluo}
%{L. Luo} and D.E. Miller,
%``Near Optimal LQR Performance for Uncertain First Order Systems'',
%{\em International Journal of Adaptive Control and Signal Processing},
%pp. 349 -- 368, May 2004.

%\bibitem{mart}
%B. M{\aa}rtensson,
%{The Order of any Stabilizing Regulator is Sufficient
%a priori Information for Adaptive Stabilization},
%{\em Systems and Control Letters}, pp. 87 -- 91, 1985.

\bibitem{rick2}
R.H. Middleton et. al,
``Design Issues in Adaptive Control",
{ \em IEEE Transaction on  Automatic Control},
  vol. 33, no. 1, pp. 50--58, 1988.

\bibitem{rick}
R.H. Middleton and G.C. Goodwin,
  ``Adaptive Control of Time-Varying Linear Systems'',
  { \em IEEE Transaction on  Automatic Control},
   vol. 33, no. 2, pp. 150--155, 1988.

\bibitem{miller03}
D.E. Miller,
``A New Approach to Model Reference Adaptive Control'',
{\em IEEE Transactions on Automatic Control},
{AC-48}, pp. 743-757, May 2003.

\bibitem{miller06}
D.E. Miller,
``Near Optimal LQR Performance for a Compact Set of Plants'',
{\em IEEE Trans. on Aut. Cont.},  pp. 1423 - 1439,
Sept. 2006.

\bibitem{scl}
D.E. Miller,
``A Parameter Adaptive
Controller Which Provides Exponential Stability: The First Order Case",
{\em Systems  and Control Letters}, vol. 103, pp. 23 -- 31, May 2017.


\bibitem{CCTA}
D.E. Miller, ``Classical Discrete-Time Adaptive Control Revisited:
Exponential Stabilization",
submitted to the {\em 1st IEEE Conference on Control Technology and Applications},August 2017.
This is posted at arXiv:1705.01494

%\bibitem{tongwen}
%D.E. Miller and T. Chen,
%``Simultaneous Stabilization with Near Optimal H-{inf} Performance'',
%{\em IEEE Transactions on Automatic Control}, {AC-47}, pp. 1986 -- 1998,
%Dec. 2002.

%\bibitem{miller91}
%D.E. Miller and E.J. Davison,
%``An Adaptive Controller which Provides an Arbitrarily
%Good Transient and Steady-State Response'',
%{\em IEEE Transactions on Automatic Control}, {\bf AC-36}, pp. 68 -- 81, 1991.

\bibitem{miller89}
D.E. Miller and E.J. Davison,
``An adaptive controller which provides Lyapunov stability",
{\em IEEE Transactions on Automatic Control}, {\bf AC-34}, pp. 599 -- 609, 1989.

%\bibitem{naghmeh}
%D.E. Miller and {\em N. Mansouri},
%``Model Reference Adaptive Control Using
%Simultaneous Probing, Estimation, and Control'',
%{\em IEEE Transactions on Automatic Control}, Vol. 55, No. 9,
%pp. 2014 -- 2029, 2010.

%\bibitem{mauro}
%D.E. Miller and {M. Rossi},
%``Simultaneous Stabilization with Near Optimal LQR Performance",
%{\em IEEE Transactions on Automatic Control}, {AC-46}, pp. 1543 --
%1555, Oct. 2001.

%\bibitem{qiuli}
%D.E. Miller and L. Qiu,
%``Near Optimal $H_{\infty}$ Tracking for a Compact Set of Plants'',
%{\em Proceedings of the IEEE 2005 Conference on Decision and Control},
%Seville, Spain, Dec. 2005.

%\bibitem{julierw}
%D.E. Miller and {\em J.R. Vale},
%``Pole Placement Adaptive Control With Persistent Jumps in
%the Plant Parameters'', {\em Mathematics of Control, Signals, and
%Systems}, vol. 26, Issue 2, pp. 177 -- 214, June, 2014.




\bibitem{morse1980}
A.S. Morse,
  ``Global Stability of Parameter-Adaptive Control Systems'',
  {\em IEEE Transactions on Automatic Control}, vol. AC-25, pp. 433--439, 1980.

\bibitem{morse96}
A.S. Morse,
``Supervisory Control of Families of Linear Set-Point Controllers -
Part 1: Exact Matching'',
{\em IEEE Transactions on Automatic Control},  {\bf AC-41},
pp. 1413 -- 1431, 1996.

\bibitem{morse97}
A.S. Morse,
``Supervisory Control of Families of Linear Set-Point Controllers -
Part 2: Robustness'',
{\em IEEE Transactions on Automatic Control},  {\bf AC-42},
pp. 1500 -- 1515, 1997.



\bibitem{naik}
S.M. Naik, P.R. Kumar and B.E. Ydstie,
``Robust Continuous-Time Adaptive Control by Parameter Projection",
{\em IEEE Transactions on Automatic Control}, vol. AC-37, pp. 182--297, 1992.

\bibitem{Nar87}
K.S. Narendra and A.M. Annawswamy,
``A New Adaptive Law for Robust Adaptation Without Persistent Excitation",
{\em IEEE Transactions on Automatic Control}, vol. AC-32, pp. 134 --145, 1987.


\bibitem{narendra}
K.S. Narendra and J. Balakrishnan,
``Adaptive Control Using Multiple Models'',
{\em IEEE Transactions on Automatic Control},  {\bf AC-42},
pp. 171 -- 187, 1997.




\bibitem{Narendra1980}
K.S. Narendra and Y.H. Lin,
  ``Stable Discrete Adaptive Control'',
  {\em IEEE Transactions on Automatic Control}, vol. AC-25, no. 3, pp. 456--461, 1980.

\bibitem{Narendra1980_pt2}
K.S. Narendra,  Y.H. Lin, and L.S. Valavani,
  ``Stable Adaptive Controller Design, Part II:
Proof of Stability'',
  {\em IEEE Transactions on Automatic Control},
  vol. AC-25, pp. 440--448, 1980.

%  Removed since it focusses on persient excitation.
%\bibitem{narendra86}
%K. Narendra  and  A. Annaswamy,
%  ``Robust adaptive control in the presence of bounded disturbances'',
%  {\em IEEE Transactions on Automatic Control},
%  vol. 31,  pp. 306--315, 1986.

\bibitem{rohrs}
C.E. Rohrs et al., ``Robustness of Continuous-Time Adaptive Control
Algorithms in the
Presence of Unmodelled Dynamics'',
{\em IEEE Transactions on Automatic Control},  {AC-30},
pp. 881 -- 889, 1985.

\bibitem{Tsakalis4}  K.S. Tsakalis and P.A. Ioannou,   ``Adaptive Control of Linear Time-Varying Plants: A New Model Reference Controller Structure'',  {\em IEEE Transaction on  Automatic Control}, vol. 34, no. 10, pp. 1038--1046, 1989.

\bibitem{julie11}
{\em J.R. Vale} and D.E. Miller,
``Step Tracking in the Presence of Persistent Plant Changes'',
{\em IEEE Transactions on Automatic Control},
pp. 43 -- 58, January, 2011.

\bibitem{vu2}
L. Vu, D. Chatterjee, and D. Liberzon,
``Input-to-State Stability of Switched Systems and Switching Adaptive Control'',
{\em Automatica}, vol. 43, pp. 639 -- 646, 2007.

\bibitem{vu}
L. Vu and D. Liberzon,
``Switching Adaptive Control of Time-Varying Plants'',
{\em IEEE Trans. on  Automat. Control}, pp. 27 -- 42,
Jan. 2011.

\bibitem{wen}
C. Wen, ``A Robust Adaptive Controller With Minimal Modifications for Discrete Time-Varying
Systems", {\em IEEE Transactions on Automatic Control}, vol. AC-39, No. 5, pp. 987 --991,
1994.

\bibitem{wenhill}
C. Wen and D.J. Hill,
``Global Boundedness of Discrete-Time Adaptive Control Using
Parameter Projection", {\em Automatica}, vol. 28, No. 2, pp. 1143 --1158, 1992.

%\bibitem{willems}
%J.C. Willems and C.I. Byrnes,
%``Global Adaptive Stabilization" in the Absence of Information
%on the Sign of the High Frequency Gain,
%{\em Lecture Notes in Control and Inf. Sciences 62},
%Springer -- Verlag: New York, pp. 49 -- 57, 1984.

\bibitem{ydstie}
B.E. Ydstie, ``Stability of Discrete-Time MRAC Revisited'',
{\em Systems and Control Letters}, vol. 13, pp. 429-439, 1989.

\bibitem{ydstie2}
B.E. Ydstie, ``Transient Performance and Robustness of Direct Adaptive Control'',
{\em IEEE Trans. on  Automat. Control}, vol. 37, No. 8, pp. 1091 -- 1105, 1992.


\bibitem{rick3}
P.V. Zhivoglyadov, R.H. Middleton, and M. Fu,
``Further Results on Localization-Based Switching Adaptive Control'',
{\em Automatica}, vol. 37, pp. 257 -- 263, 2001.

\bibitem{zhou_book}
K. Zhou, J.C. Doyle and K. Glover,
{\em Robust and Optimal Control}, Prentice Hall, New Jersey, 1995.


\end{thebibliography}
\end{document}